\documentclass[11pt]{preprint}
\usepackage[full]{textcomp}
\usepackage{pgfplots}
\pgfplotsset{compat=1.17}
\usepgfplotslibrary{groupplots}
\usepackage[osf]{newtxtext} 
\usepackage[cal=boondoxo]{mathalfa}
\usepackage{colortbl}
\usepackage[normalem]{ulem}
\usepackage{amsmath}
\usepackage{pgf}
\usepackage{xcolor}
\usepackage{comment}

\usepackage{amssymb}
\usepackage{mathtools}
\usepackage{hyperref}
\usepackage{breakurl}
\usepackage{mhenvs}
\usepackage{mhequ} 
\usepackage{mhsymb}
\usepackage{booktabs}
\usepackage{tikz}
\usepackage{tcolorbox}
\usepackage{mathrsfs}
\usepackage[utf8]{inputenc}
\usepackage{longtable}
\usepackage{wrapfig}

\usepackage{microtype}
\usepackage{comment}
\usepackage{wasysym}
\usepackage{centernot}
\usepackage{enumitem}
\usepackage{bm}
\usepackage{stackrel}
\usepackage{graphicx}
\usepackage{subcaption}

\makeatletter
\newcommand{\globalcolor}[1]{%
	\color{#1}\global\let\default@color\current@color
}
\makeatother

\usetikzlibrary{calc}
\usetikzlibrary{decorations}
\usetikzlibrary{positioning}
\usetikzlibrary{shapes}
\usetikzlibrary{external}

\newif\ifdark
\darkfalse

\ifdark

\definecolor{darkred}{rgb}{0.9,0.2,0.2}
\definecolor{darkblue}{rgb}{0.7,0.3,1}
\definecolor{darkgreen}{rgb}{0.1,0.9,0.1}
\definecolor{franck}{rgb}{0,0.8,1}
\definecolor{pagebackground}{rgb}{.15,.21,.18}
\definecolor{pageforeground}{rgb}{.84,.84,.85}
\pagecolor{pagebackground}
\AtBeginDocument{\globalcolor{pageforeground}}
\definecolor{symbols}{rgb}{0,.7,1}
\colorlet{connection}{red!80!black}
\colorlet{boxcolor}{blue!50}

\else

\definecolor{darkred}{rgb}{0.7,0.1,0.1}
\definecolor{darkblue}{rgb}{0.4,0.1,0.8}
\definecolor{darkgreen}{rgb}{0.1,0.7,0.1}
\definecolor{franck}{rgb}{0,0,1}
\definecolor{pagebackground}{rgb}{1,1,1}
\definecolor{pageforeground}{rgb}{0,0,0}
\colorlet{symbols}{blue!90!black}
\colorlet{connection}{red!30!black}
\colorlet{boxcolor}{blue!50!black}

\fi

\def\slash{\leavevmode\unskip\kern0.18em/\penalty\exhyphenpenalty\kern0.18em}
\def\dash{\leavevmode\unskip\kern0.18em--\penalty\exhyphenpenalty\kern0.18em}

\DeclareMathAlphabet{\mathbbm}{U}{bbm}{m}{n}

\DeclareFontFamily{U}{BOONDOX-calo}{\skewchar\font=45 }
\DeclareFontShape{U}{BOONDOX-calo}{m}{n}{
	<-> s*[1.05] BOONDOX-r-calo}{}
\DeclareFontShape{U}{BOONDOX-calo}{b}{n}{
	<-> s*[1.05] BOONDOX-b-calo}{}
\DeclareMathAlphabet{\mcb}{U}{BOONDOX-calo}{m}{n}
\SetMathAlphabet{\mcb}{bold}{U}{BOONDOX-calo}{b}{n}

\setlist{noitemsep,topsep=4pt,leftmargin=1.5em}

\DeclareMathAlphabet{\mathbbm}{U}{bbm}{m}{n}

\DeclareMathAlphabet{\mcb}{U}{BOONDOX-calo}{m}{n}
\SetMathAlphabet{\mcb}{bold}{U}{BOONDOX-calo}{b}{n}
\DeclareFontFamily{U}{mathx}{\hyphenchar\font45}
\DeclareFontShape{U}{mathx}{m}{n}{
	<5> <6> <7> <8> <9> <10>
	<10.95> <12> <14.4> <17.28> <20.74> <24.88>
	mathx10
}{}
\DeclareSymbolFont{mathx}{U}{mathx}{m}{n}
\DeclareMathSymbol{\bigtimes}{1}{mathx}{"91}

\providecommand{\figures}{false}
{ \ifthenelse{\equal{\figures}{false}} {#1}{\[ {\rm Figure \ missing !} \]} }{}

\newtheorem{example}{Example}

\newtheorem{scheme}{Scheme}[section]

\newtheorem{hypothesis}{Hypothesis}[section]

\definecolor{oxfordblue}{rgb}{0.0, 0.13, 0.28}
\tikzstyle{leaf}=[circle, draw=black, fill=gray!20, inner sep = .1pt, font=\tiny, minimum size=3mm]
\tikzstyle{inner}=[circle, minimum size=1.5mm, fill=oxfordblue, inner sep =0]
\tikzstyle{conj} = [thick, color=maroon, dotted,dash pattern=on 1pt off 1pt]
\tikzstyle{conj-int} = [color=cyan, dotted]
\tikzstyle{int} = [color=cyan]
\tikzstyle{positive} = [thick, color=maroon]
\tikzstyle{Phi} = [thick, color=maroon, snake=snake,
segment length=5pt,
segment amplitude=1pt]
\tikzstyle{Phi-conj} = [thick, color=maroon, snake=snake,
segment length=5pt,
segment amplitude=1pt,dashed]
\tikzstyle{stoch-int} = [color=cyan, snake=zigzag,
segment amplitude=.2mm,
segment length=.5mm]
\tikzstyle{stoch-conj-int} = [color=cyan, snake=triangles,
segment amplitude=.25mm,
segment length=.8mm]
\definecolor{modebeige}{rgb}{0.59, 0.44, 0.09}
\definecolor{maroon}{rgb}{0.5, 0.0, 0.0}
\definecolor{mayablue}{rgb}{0.45, 0.76, 0.98}

\def\mcO{\mathcal{O}}
\def\mbbR{\mathbb{R}}

\tikzstyle{tinydots}=[dash pattern=on \pgflinewidth off \pgflinewidth]
\tikzstyle{superdense}=[dash pattern=on 4pt off 1pt]

\newcommand{\mcM}{\mathcal{M}}
\newcommand{\mcE}{\mathcal{E}}

\newcommand{\mcA}{\mathcal{A}}

\newcommand{\mcL}{\mathcal{L}}

\newcommand{\mcF}{\mathcal{F}}

\newcommand{\mcK}{\mathcal{K}}

\newcommand{\mcP}{\mathcal{P}}

\newcommand{\mbbE}{\mathbb{E}}

\newcommand{\mbbN}{\mathbb{N}}
\newcommand{\mbbP}{\mathbb{P}}

\newcommand{\mbbT}{\mathbb{T}}

\newcommand{\mbbZ}{\mathbb{Z}}

\newcommand{\sbg}{\boldsymbol{g}}
\newcommand{\sbf}{\boldsymbol{f}}

\newcommand{\bmu}{\boldsymbol{\mu}}
\newcommand{\bnu}{\boldsymbol{\nu}}
\newcommand{\bK}{\boldsymbol{K}}
\newcommand{\bL}{\boldsymbol{L}}

\newcommand{\sbv}{\boldsymbol{v}}
\newcommand{\sbw}{\boldsymbol{w}}

\newcommand{\del}{\delta}

\newcommand{\mfF}{\mathfrak{F}}

\newcommand{\mfP}{\mathfrak{P}}

\newcommand{\mfR}{\mathfrak{R}}

\newcommand{\mfS}{\mathfrak{S}}

\def\${|\!|\!|}

\def\inner#1{{\left<#1\right>}}

\newcommand{\pa}{\partial}
\newcommand{\rbr}[1]{\left(#1\right)}

\newcommand{\sqbr}[1]{\left[#1\right]}
\newcommand{\abs}[1]{\left| #1 \right|}
\newcommand{\norm}[1]{\left| \left| #1 \right| \right|}

\newcommand{\set}[1]{\left\lbrace  #1 \right\rbrace }

\newcommand{\ti}[1]{\tilde{#1}}

\newcommand{\lam}{\lambda}

\newcommand{\grad}{\nabla}

\newcommand{\lap}{\Delta}

\newcommand{\pint}[1]{\the\numexpr #1 \relax}

\newenvironment{DIFnomarkup}{}{} 

\theorembodyfont{\rmfamily}

\newfont{\indic}{bbmss12}

\def\Nabla_#1{\nabla_{\!#1}}

\let\eps\varepsilon

\def\eqref#1{(\ref{#1})}

\makeatletter 
\newcommand*{\bigcdot}{}
\DeclareRobustCommand*{\bigcdot}{%
	\mathbin{\mathpalette\bigcdot@{}}%
}
\newcommand*{\bigcdot@scalefactor}{.5}
\newcommand*{\bigcdot@widthfactor}{1.15}
\newcommand*{\bigcdot@}[2]{%
	\sbox0{$#1\vcenter{}$}
	\sbox2{$#1\cdot\m@th$}%
	\hbox to \bigcdot@widthfactor\wd2{%
		\hfil
		\raise\ht0\hbox{%
			\scalebox{\bigcdot@scalefactor}{%
				\lower\ht0\hbox{$#1\bullet\m@th$}%
			}%
		}%
		\hfil
	}%
}
\makeatother

\tcbset
{colframe=boxcolor,colback=symbols!7!pagebackground,coltext=pageforeground,
	fonttitle=\bfseries,nobeforeafter,center title,size=fbox,boxsep=1.5pt,
	top=0mm,bottom=0mm,boxsep=0mm,tcbox raise base}

\def\two{{\<generic>\kern0.05em\<genericb>}}
\def\twoI{{\<Ito>\kern0.05em\<Itob>}}

\def\mail#1{\burlalt{#1}{mailto:#1}}

\begin{document}
	\usetikzlibrary{snakes}
	\title{Low regularity symplectic schemes for stochastic NLS}\author{Jacob Armstrong-Goodall$^1$, Yvain Bruned$^2$}
	\institute {Maxwell Institute for Mathematical Sciences, University of Edinburgh \and IECL (UMR 7502), Université de Lorraine
		\\
		Email:\ \begin{minipage}[t]{\linewidth}
			\mail{j.a.armstrong-goodall@sms.ed.ac.uk} \\ \mail{yvain.bruned@univ-lorraine.fr}.
	\end{minipage}}
	
	\maketitle
	
	\begin{abstract}
		We introduce a class of symplectic resonance based schemes for Schr\"odinger's equation in dimension one, building on the work in \cite{AB2023} wherein resonance based numerical schemes were developed in the context of dispersive PDE driven by time dependent, or space-time dependent, coloured noise. We work primarily with a cubic nonlinearity, advancing the approach introduced in \cite{M.S2023} for deriving symplectic schemes in the deterministic setting. As an example of such a scheme we derive the resonance based midpoint rule for the Stochastic NLS and analyse its convergence properties.
	\end{abstract}

	\setcounter{tocdepth}{2}
	\tableofcontents
	
	\section{Introduction}
	
	Symplectic numerical schemes are those that preserve certain geometrical structures inherent in the solution of a differential equation. Such integrators can be seen as a discrete process possessing the same qualitative properties as the associated Hamiltonian system. The geometrical structures in question often arise in physics representing physical quantities, such as mass or energy, that are preserved over the time evolution of a system. Symplecticity is closely related to integrability which is the study of preserved properties solutions to differential equations. Probably the most famous example of an integrable system from the PDE literature is the KdV equation. This equation preserves mass, momentum, energy and many other quantities. In this work we will study numerical schemes for the nonlinear Schr\"odinger equation driven by space-time white noise,
	\begin{equation*}
		i\pa_t u + \pa_x^2 u = \lam\abs{u}^{2}u  + \kappa u\circ \Phi\xi(t,x),
	\end{equation*}
	with solutions on the one dimensional torus $ \mbbT $. The schemes we introduce will preserve two important quantities, the first being the mass
	\begin{equation}
		\int_\mbbT |u(s,x)|^2 dx
	\end{equation}
	and the second being the stochastic perturbed Hamiltonian system,
	\begin{align}
		\dot{u} = -\frac{\partial H_0}{\partial \bar{u}} - \frac{\partial H_1}{\partial \bar{u}} \circ \Phi \xi \quad
		\dot{\bar{u}} = \frac{\partial H_0}{\partial u} + \frac{\partial H_1}{\partial u} \circ \Phi \xi
	\end{align}
	where
	\begin{equ}
		H_0(u) =  \frac12\int\abs{\grad u}^2 dx + \frac{\lam}{4}\int \abs{u}^{4} dx, \quad
		H_1(u) = \frac{\kappa}{2}\int \abs{u}^2 dx.
	\end{equ}
	This can formally be thought of as the Hamiltonian of the deterministic system $ H_0(u) = \frac12\norm{\grad u}_{L^2} + \frac{\lam}{4}\norm{u}_{L^4}  $, with the addition of the work done by the noise, $ H_1\circ \Phi \xi $. 
	
	The study of numerical methods for Hamiltonian systems is a well developed field making use of various techniques including variational, symplectic and multi-symplectic methods, see \cite{B.R2006} for a robust review article on the topic. More recently resonance based schemes for dispersive PDE have been introduced to deal with the low regularity regime, see \cite{B.S2022} for high order local error analysis of a very general class of equations. In \cite{B.B.S2022}, more general non-linearities are considered by avoiding use of the Fourier transform, and instead the analysis is performed via nested commutators. Symmetric-in-time resonance based schemes are introduced in \cite{B.B.M.S2023} and play a pivotal role in the construction of symplectic schemes, while in \cite{M.S2023}, Runge-Kutta methods are derived for KdV and Schr\"odinger's equation and the error analysis of the midpoint rule is performed as an example.
	
	In the case of Schr\"odinger's equation the idea in \cite{M.S2023} is to first use polynomial interpolation to give a low regularity, symplectic, discretisation of the kernel, leading to a map of the intial conditions to an approximation of the solution at a given time. Secondly, the S-stage Runge-Kutta method is applied to the Hamiltonian system describing the evolution of the equation. Combining these concepts gives rise to a class of implicit Runge-Kutta schemes which at low orders have desirable regularity properties. These are the first instance of symplectic schemes designed to deal with low regularity initial data. The resonance based schemes introduced in \cite{B.S2022} were applied in the stochastic setting in \cite{AB2023} where strong and pathwise global error analysis was performed for schemes approximating solutions to equations of the form, 
	\begin{equation}
		\begin{aligned}
			& i \pa_t u(t,x) +   \mathcal{L}\left(\nabla \right) u(t,x) =\vert \nabla\vert^\alpha p\left(u(t,x), \overline u(t,x)\right) \\ & +\vert \nabla\vert^\beta  f(u(t,x),\overline u(t,x)) 	\Phi \xi(t,x), \quad
			u(0,x) = v(x).
		\end{aligned} 
	\end{equation}
	where $\Phi \xi(t,x) $ is a white noise in time and coloured in space. In the sequel, we will also use the notation $ \Phi dW(t,x) $.  We find that the fixed point convergence, and stability of the scheme both require a restrictive bound on the time-step which is random and proportionate to the norm of the initial conditions. Due to this obstacle we only study the local error of the scheme.

	In \cite{AB2023} it is shown that one can obtain an order $ \sqrt{t} $ approximation to the SNLS whith initial data $ v\in H^1 $ and $ \Phi \xi(t,x)\in H^2 $ in space. At higher order it is shown that no improvement is offered over traditional schemes due to lack of path-wise integrability of the stochastic integrals arising in the iterations of Duhamel's formula. In this work we note that symplecticity in the stochastic setting is only satisfied by the averaged energy. Perhaps for this reason, the study of stochastic Hamiltonian systems began only recently, for instance, schemes derived to preserve symplectic structures were proposed by Milstein, Repin and Tetryakov in \cite{MRT2002} for SDE, these authors overcame the restriction on the time-step by using a truncated noise, this however does not lead to a strongly convergent scheme. In \cite{D.D2002}, Debussche and Di Menza performed numerical experiments to explore the impact of the spatial covariance of the noise on the blow-up of the solutions to the focusing cubic SNLSE. For this they employ the stochastic midpoint, the low regularity version of which is derived in this paper. Debussche and de Bouard then study the discretisation used in these experiments in \cite{B.D2004}, wherein they observe, as we do, that the smallness condition for stability is random and too restrictive. Interestingly, with regards to the simulations in  \cite{D.D2002}, they state that the restriction caused no apparent problems in the implementation of the scheme and we reproduce this result our implementation of the schemes. The same authors also prove the convergence in probability of a mass preserving scheme for SPDE with Itô type noise in \cite{B.D2006}. Symplectic Runge-Kutta and mean square convergence analysis of midpoint rule for SNLSE is performed in \cite{CheHo2016} and in \cite{H.W2019} where both symplectic and multi-symplectic methods are studied. Strong convergence for the SNLSE with random initial data is studied in \cite{CH2018, CuHoL2017} where space-time discretisation is considered for the Stratanovich product in the former, while in the latter the time discretisation is studied for the Itô noise and damping proportional to the covariance operator. Strong convergence of an operator splitting method for the cubic SNLS with space-time Stratonovich noise (coloured in space, white in time), on $ \mbbR $  is performed in \cite{CHLZ2019}. The works focusing on obtaining strong convergence results rely on exponential integrability conditions to perform the stability and error analysis. The exponential integrability of not only the scheme, but the solution itself, seems to be a key obstacle in obtaining results about strong convergence, especially in dimensions greater than one. When it comes to implicit schemes however the bound on the time-step is restrictive enough to imply the exponential integrability but is so restrictive that stability can never be ensured.

	The overarching goal in this work is to derive a class of symplectic resonance based Runge-Kutta schemes for Schr\"odingers equation in one dimension. In contrast to \cite{AB2023} we will work with equations with noise of the Stratonovich type. This is because the symplectic discretisation of the nonlinear term is necessarily implicit and Itô type noises are inherently linked to explicit discretisations of the noise. Another difference is that for the deterministic integral arising in Duhamel's formula, instead of straightforward Taylor expansion of the operator, we will employ techniques used in \cite{M.S2023} via the symplectic low regularity approximation to the kernel 
	\begin{equs} \label{kernel_symplectic}
		\mcK(s;k,k_1,k_2,k_3):=e^{-2iskk_1} + e^{2isk_2k_3} - 1.
	\end{equs}
	For the stochastic convolution we will use a similar approximation but note that as observed in \cite{AB2023}, the integrals can not be studied in the pathwise sense, meaning that strongly convergent schemes can not be constructed this way. We will also employ novel iterations of Duhamel's formula, first introduced in \cite{M.S2023} as follows:\[u(t) = e^{i(t-s)\lap}u(s) - ie^{it\lap}\int_s^t F(u(r),\bar{u}(r))dr - ie^{it\lap}\int_s^tu\Phi\circ dW(x,r)\]which, when used in conjunction with  the standard mild form to derive a midpoint estimate to $ u(t) $, produces cancellations in the second iteration of the stochastic convolution. This implies the second iteration of the stochastic convolution in Duhamel's formula is zero up to a term of order $ t^2 $, which allows us to construct low order schemes without this term. The above considerations lead to the following class of schemes for the SNLSE
	\begin{align*}
		u_{n+1} &= e^{it \pa_x^2} u_n + t e^{it \pa_x^2} \sum_{\alpha \in D} b^{(0)}_{\alpha} K_{\alpha} + \sqrt{t} e^{it \pa_x^2} \sum_{\alpha \in D} b^{(1)}_{\alpha} L_\alpha \\
		U_{\alpha}(u_n) &= u_n + t\sum_{\ti{\alpha} \in D} a^{(0)}_{\alpha,\ti{\alpha}} K_{\ti{\alpha}} + \sqrt{t}\sum_{\ti{\alpha}\in D} a^{(1)}_{\alpha,\ti{\alpha}} L_{\ti{\alpha}}
	\end{align*}
	where 	$\alpha \in D = \set{(p,q,r):0 \leq p,q,r \leq S}$ and 
	\begin{equation*}
		K_{\alpha} = \mathfrak{F}^{2d}_{p}\rbr{t;c_{q};U_{\alpha}(u_n)},\quad
		L_{\alpha} = \mathfrak{P}^{2d}_{p}\rbr{t;c_{q};U_{\alpha}(u_n)}.\
	\end{equation*}
	Here $ K_\alpha $ and $ L_\alpha $ represent the discretisation of the deterministic and stochastic integrals in Duhamel's formulation and $ \alpha = (p,q,r) $ is the coefficient arising in the S-stage Runge-Kutta iteration of the Hamiltonian. We will prove that the mass, and the sample average of the energy, are preserved under the correct assumptions on the noise and the coefficients of the discretisation. As an example of such a scheme we will derive the stochastic resonance midpoint rule 
	\begin{equation*}
		S^{n+1}_t \rbr{u_0} = u_{n+1} = e^{it\pa_x^2}u_n + e^{it\pa_x^2}\rbr{t\mathfrak{F}_{0}\rbr{t;1;u_{n+\frac12}} + \sqrt{t}\mathfrak{P}_{0}\rbr{t;1;u_{n+\frac12}}}
	\end{equation*}
	where
	\begin{equation*}
		u_{n+\frac12} = \frac{u_n + e^{-it\pa_x^2}u_{n+1}}{2},
	\end{equation*}
	which is the resonance based version of the scheme in \cite{CheHo2016,D.D2002}. We prove the strong and pathwise convergence of the midpoint scheme and note  that the regularity requirements are improved as a result of the careful treatment of the oscillatory integrals. It should be pointed out that these schemes can be applied to other dispersive equations such as the wave equation, Manakov system or KdV by deriving the corresponding dicretisations of the terms in Duhamel's formula.  The only limitation is to be able to perform the symplectic discretisation of a kernel of the type \eqref{kernel_symplectic}.
		Let us mention that there were recent developments in the error analysis in the deterministic setting (see \cite{RouSc2022, LiWu2025, CaoYa2025}) where the authors were looking at rough initial data with regularity below $ H^1 $. One can also play the same  game for random initial data and one expects to go down to lower regularity as it is the case for the well-posedness theory.

	Let us briefly summarise our main results.  We present a new class of stochastic Runge-Kutta schemes \eqref{scheme:stoch-runge-kutta} for the one-dimensional SNLSE. One specific example, Scheme~\eqref{ex:stoch-res-midpoint}, which employs the midpoint iteration and a symplectic approximation of the kernel \eqref{kernel_symplectic}, exhibits low regularity. In Theorem \ref{thm:Symplecticity}, by imposing Hypothesis \ref{hyp:coefficient}, we can show that the previously introduced schemes are symplectic. In Section \ref{sec::6}, the analysis of the existence, uniqueness, and local error is performed on the low regularity scheme \eqref{ex:stoch-res-midpoint}, but we observe that its implicit nature imposes strict conditions on the step size, which prevents stability in the presence of noise. 
	
	Finally, let us outline the content of each section. In Section \ref{sec::2}, the derivation of Runge-Kutta schemes for deterministic PDE is summarised from \cite{M.S2023}, see equation \eqref{eq:scheme_nlse} and Theorem \ref{condition_deterministic}. This idea is then extended to the stochastic setting in Section \ref{sec::4} where the general class of implicit Runge-Kutta schemes is presented in Scheme \ref{scheme:stoch-runge-kutta}. We use techniques introduced in \cite{B.B.M.S2023} to make the key observation that the second iteration of the stochastic integral disappears up to order $ t^2 $ and introduce hypotheses \ref{hyp:Phi} and \ref{hyp:coefficient}, required for the conservation of energy and mass. Finally, in this section, the midpoint rule is derived explicitly in example \ref{ex:stoch-res-midpoint}. 
	
	In Section \ref{sec::5}, we prove that the class of schemes derived in Section \ref{sec::4} conserves the quantities introduced in Section \ref{sec::Hamiltonian-Form}. This result is summarised in Theorem \ref{thm:Symplecticity}, which is then proved by combining the technical result Lemma \ref{lemma:vwedgeGiszero} and the conservation laws for mass, Proposition \ref{prop:conservation-of-mass}, and energy, Proposition \ref{prop:conservation-of-energy}. The calculations are similar to those of \cite{H.W2019} and \cite{M.S2023}. We note that in the stochastic setting, symplecticity occurs only after taking into account the averaging of the Fourier coefficients.
	
	In Section \ref{sec::6}, we first study the stability of the scheme introduced in Example \eqref{ex:stoch-res-midpoint}, the result of which is presented in Theorem \ref{thm:stability}. Next, we study the convergence of the fixed point iterations, which is the result expressed in Theorem \ref{thm:fixedpoint-convergence}. We observe that these properties rely on the condition that the time step is sufficiently small, but this cannot be ensured due to the presence of randomness. For this reason, we cannot prove the convergence of the scheme. We find that we can still prove the local error is of order $ t^{3/2} $ using the cancellations that occur in the midpoint iterations of the Duhamel formulation; this is the statement of Theorem \ref{prop:local-error}. 
	\subsection*{Acknowledgements}
	
	{\small
		We would like to thank Jianbo Cui who spot the importance of exponential integrability conditions for getting strong convergence results.	Y. B. gratefully acknowledges funding support from the European Research Council (ERC) through the ERC Starting Grant Low Regularity Dynamics via Decorated Trees (LoRDeT), grant agreement No.\ 101075208. Views and opinions expressed are however those of the author(s) only and do not necessarily reflect those of the European Union or the European Research Council. Neither the European Union nor the granting authority can be held responsible for them. 
		J. A. G. is supported by the EPSRC Centre for Doctoral Training in Mathematical Modelling, Analysis and Computation (MAC-MIGS) funded by the UK Engineering and Physical Sciences Research Council (grant EP/S023291/1), Heriot-Watt University and the University of Edinburgh.
	}

	\section{Hamiltonian Equations}\label{sec::Hamiltonian-Form}
	Our focus will primarily be on the stochastic cubic Schr\"odinger equation on the one dimensional torus $ \mbbT. $
	\paragraph{Stochastic NLS Equation (SNLSE):}
	\begin{equation}\label{eq::snlse}
		i\partial_t u(t,x) + \pa_x^2 u(t,x) - \lam\abs{u(t,x)}^{2}u(t,x) = \kappa u(t,x)\circ\Phi \xi(t,x),\quad u(0,x)=v(x).
	\end{equation}	
	Here $ \Phi:L^2(\mbbT) \rightarrow  H^s(\mbbT) $ is a Hilbert-Schmidt operator. We denote $ \xi(t,x) = \frac{\pa^2 W(t,x)}{\pa t\pa x} $ and impose that it satisfies the identity,
	\begin{equation*}
		\mbbE\sqbr{\xi(t,x) \xi(s,y)} = \delta(t-s)\delta(x-y),
	\end{equation*}
	which is the definition of space-time white noise. The operator $ \Phi\Phi^* $ can be interpreted as the spatial covariance of $\xi(t,x)$, specifically, if $ \Phi $ is defined through a kernel $ K $ as follows:
	\begin{equation*}
		\Phi u(x) = \int_{\mbbT}K(x,y)u(y)dy.
	\end{equation*}
	Then the correlation of the noise is given by
	\begin{equation*}
		\mbbE[\Phi\xi(t,x)\Phi\xi(s,y)] = c(x,y)\delta(t-s)
	\end{equation*}
	with
	\begin{equation*}
		c(x,y) = \int_{\mbbT^d}K(x,z)K(y,z)dz,
	\end{equation*}
	reflecting the fact that increasing the spatial covariance is equivalent to regularisation. 
	
	The mild formulation or the Duhamel formulation of the cubic SNLSE is 
	\begin{equation} \label{eq::intro_SNLS_duhamel} 
		\begin{aligned}
			u(t) &= e^{i t \Delta}v - i\lam e^{i t \Delta}\int_0^t e^{-i s \Delta}u(s)\abs{u(s)}^2ds \\ &- i\kappa e^{i t \Delta}\int_0^t e^{-i  s \Delta} u(s)\circ\Phi dW(s).
		\end{aligned}
	\end{equation} 
	Consider the infinite-dimensional phase space $ M $ with coordinates $ \mu_k $ and momenta $ \bar{u}_k $, for $ k\in\mbbZ$. Then the Hamiltonian is a function $ H:M\rightarrow \mbbR $ and the solution to the Hamiltonian system is a curve $ \rbr{u_k(t),\bar{u}_k(t)} $ in $ M $ which obeys the Hamiltonian equations of motion:
	\begin{equation}\label{eq:geo-struc-det}
		\dot{u}_k = \frac{\pa H}{\pa \bar{u}_k},\quad \dot{\bar{u}}_k = -\frac{\pa H}{\pa u_k}.
	\end{equation}  
	When the Hamiltonian is time independent, $ \pa H/\pa t = 0  $, the value of the function $ H(u,\bar{u}) $ is conserved over the evolution of the solutions of the equation. This conservation law can also be expressed via the symplectic form,
	\begin{equation}\label{eq:symplectic-form}
		\omega(u,\bar{u})=\int du(t, x)\wedge d\bar{u}(t, x)dx = \int dv(x)\wedge d\bar{v}(x) dx,\quad\forall t\geq0.
	\end{equation}
	The pair $ (M,\omega) $ is a symplectic manifold. This identity corresponds to the conservation of energy in physical systems.
	
	\subsubsection{The Deterministic Setting}
	First we consider the deterministic setting in which $ \kappa = 0 $,
	\begin{equation}\label{eq::nlse}
		i\partial_t u(t,x) + \pa_x^2 u(t,x) - \lam\abs{u(t,x)}^{2}u(t,x) = 0,\quad u(0,x)=v(x).
	\end{equation}
	When considering a periodic boundary, the equation \eqref{eq::nlse} arises from the following Hamiltonian: 
	\begin{equation}\label{eq:hamiltonian}
		H(u) = \int \left(\frac12 \abs{\pa_x u}^2 + \frac{\lam}{4}\abs{u}^{4} \right) dx,
	\end{equation}
	which in physics is interpreted as the energy of the system, and in the case of the defocusing NLSE, is a conserved quantity.
	In Fourier space, one has $ \mcF(\pa_x v)(k) = i k \hat{v}_{k}$, which gives
	\begin{equation}\label{eq::intro_Ham_Fourier}
		H(u,\bar{u}) = \frac{1}{2}\sum_{k\in \mbbZ}\abs{k}^2u_k\bar{u}_k + \frac{\lam}{4}\sum_{k_1,k_2,k_3,k_4 \in \mbbZ}u_{k_1}\bar{u}_{k_{2}}u_{k_3}\bar{u}_{k_{4}}.
	\end{equation}
	The solution of \eqref{eq::nlse} corresponds to the curves $ \rbr{u_k(t),\bar{u}_k(t)} $ obeying \eqref{eq:geo-struc-det}. In the case of \eqref{eq::nlse}, the other conserved quantity is the mass,
	\begin{align}\label{eq:conserved}
		\mcM = \int \abs{u}^2 dx = \sum_{k\in\mbbZ}u_k\bar{u}_k, 
	\end{align}
	The Hamiltonian system is defined on the phase space $ \mcP_s(\mbbR) := l^2_s(\mbbR)\times l^2_s(\mbbR)$ where for $ s\in\mbbR $
	\begin{equation*}
		l^2_s(\mbbR):=\set{(a_k)_{k\in\mbbZ}\in\mbbR^\mbbZ:\sum_{k\in\mbbZ}k^{2s}\abs{a_k}^2<\infty}
	\end{equation*}
	is a Hilbert space for the standard norm $ \norm{a}_s^2=\sum_{k\in\mbbZ}\abs{k}^{2s}\abs{a_k}^2 $.
	\subsubsection{The SNLSE as a Hamiltonian System}\label{intro::SNLS-as-Ham}
	We consider generic Hamiltonian systems in one dimension of the form 
	\begin{equation}\label{eq::HamiltonainIntro}
		\dot{u}_k = \frac{\delta H}{\delta \bar{u}_k} + \frac{\delta \hat{H}}{\delta \bar{u}_k}\circ\Phi\xi,
		\quad \dot{\bar{u}}_k = -\frac{\delta H}{\delta u_k} - \frac{\delta \hat{H}}{\delta u_k}\circ\Phi\xi.
	\end{equation}
	Here $ \xi(t,x) =\dot{W}(t,x)$ is a space-time white noise and $ W(t,x) $ is a real-valued cylindrical Brownian motion defined by
	\begin{equation}\label{eq::real-cylindrical-BM}
		W(t,x) = \sum_{k\in\mbbZ}W_k(t)e^{ikx},\quad W_k(t):\Omega\rightarrow \mbbR,
	\end{equation} 
	where $ W_k(t) = \overline{W_{-k}(t)} $. The SNLSE \ref{eq::snlse} is described by \ref{eq::HamiltonainIntro} with the following identities for $ H $ and $ \hat{H} $: 
	\begin{equation}\label{eq:stoch-hamiltonian}
		H(u) =  \frac12\int\abs{\grad u}^2 dx + \frac{\lam}{4}\int \abs{u}^{4} dx, \quad
		\hat{H}(u) = \frac{\kappa}{2}\int \abs{u}^2 dx.
	\end{equation}
	As in \eqref{eq::intro_Ham_Fourier}, under the Fourier transform we see that $ H(u,\bar{u}) = H(u) $ with the stochastic part described as 
	\begin{align*}
		\hat{H}(u,\bar{u}) = \frac\kappa2\sum_{k_1\in\mbbZ}u_{k_1}\bar{u}_{k_1},\quad \Phi\xi = \Phi_{k_2} dW_{k_2}(t).
	\end{align*}
	In the stochastic setting, the mass is conserved just as in the deterministic case, but the energy is not. However, as in \cite[Proposition 2]{CheHo2016} and \cite[Proposition 4.5]{BouDe2003}, the average energy can be described as
	\begin{equation}\label{eq::averaged-energy}
		\mathbb{E}\left[H(u(t))\right] = \mathbb{E}\left[H(v)\right] + \frac{\kappa^2}{2} \sum_{\substack{k_1, k_2 \in \mathbb{Z} \\ k_1 = -k_2}} \int_0^t k_2^2 u_{k_1} \bar{u}_{k_1} \Phi_{k_2} \Phi^*_{k_2} \, ds
	\end{equation}
	and this quantity is independent of time.
	We will now describe how a symplectic form arises from the Hamiltonian system \eqref{eq:stoch-hamiltonian}.
	Given two functionals $ \mcL$ and $ \hat{H} $ in the functions $ P,Q:\mbbT\times \mbbR_+ \rightarrow \mbbR $ the generalised action integral is given by
	\begin{equation}\label{eq::action}
		\mcA = \int_{0}^T\rbr{\mcL(P,Q,\dot{P},\dot{Q}) - \hat{H}(P,Q)\circ\Phi\xi}dt.
	\end{equation}
	\begin{definition}[Functional Derivative]
		Suppose $ B $ is a Banach space and $ F $ is a functional defined on B. The functional differential of $ F $ at the point $ b\in B $ is the linear functional $ \delta F\sqbr{b,\cdot} $ in B defined such that $ \forall\phi\in B $,
		\begin{equation*}
			\delta F\sqbr{b,\phi}=\lim_{\eps\rightarrow 0}\frac{F[b+\eps\phi] - F[b]}{\eps}.
		\end{equation*}
		The functional derivative of $ F = \int_\mbbT L(b,\dot{b})dx $ is
		\begin{equation*}
			\delta F[b,\phi] = \int_\mbbT\frac{\delta F}{\delta b}(x)\phi(x)dx.
		\end{equation*} 
	\end{definition} 
	With this definition in hand, we may write the variation of the action as
	\begin{equs}\label{eq::action_variation}
		\begin{aligned}
	&	\delta\mcA = \delta\int_0^T\rbr{\mcL(P,Q,\dot{P},\dot{Q}) - \hat{H}(P,Q)\circ\Phi\xi}dt\notag \\
		&=\int_\mbbT\int_0^T\left(\frac{\delta\mcL}{\delta P}\delta P + \frac{\del\mcL}{\del Q}\del Q + \frac{\del\mcL}{\del\dot{P}}\del\dot{P} + \frac{\del\mcL}{\del\dot{Q}}\del \dot{Q} \right.
			\\ & \left. - \frac{\del \hat{H}}{\del P}\del P\circ \Phi\xi -\frac{\del \hat{H}}{\del Q}\del Q\circ \Phi\xi \right)dtdx \notag\\
		&=\int_\mbbT\int_0^T\left[ \rbr{\frac{\del\mcL}{\del P} - \frac{d}{dt}\rbr{\frac{\del\mcL}{\del \dot{P}}} - \frac{\del \hat{H}}{\del P}\circ\Phi\xi}\del P \right. \\ & \left. +\rbr{ \frac{\del\mcL}{\del Q} - \frac{d}{dt}\rbr{\frac{\del\mcL}{\del \dot{Q}} }- \frac{\del \hat{H}}{\del Q}\circ\Phi\xi}\del Q \right]dtdx.
		\end{aligned}
	\end{equs}
	The final equality follows by integration by parts, and we assume the boundary conditions
	\begin{equation*}
		\del P(T,x) = \del P(0,x) = \del Q(T,x) = \del Q(0,x) =0, \quad x\in\mbbT.
	\end{equation*}
	Hamilton's principle states that the evolution of \eqref{eq::HamiltonainIntro} over the interval $ (0,T) $ is a stationary point of the action integral, i.e., a solution to the system
	\begin{equation*}
		\del\mcA = \del\int_{0}^T\rbr{\mcL(P,Q,\dot{P},\dot{Q}) - \hat{H}(P,Q)\circ\Phi\xi}dt = 0.
	\end{equation*}
	Using \eqref{eq::action_variation}, we obtain the system
	\begin{equation}\label{eq::systemA}
		\frac{d}{dt}\rbr{\frac{\del\mcL}{\del\dot{P}}} = \frac{\del\mcL}{\del P} - \frac{\del \hat{H}}{\del P}\circ\Phi\xi,\quad \frac{d}{dt}\rbr{\frac{\del\mcL}{\del\dot{Q}}} = \frac{\del\mcL}{\del Q} - \frac{\del \hat{H}}{\del Q}\circ\Phi\xi.
	\end{equation}
	Define the functional from \eqref{eq::action} by \[\mcL(P,Q,\dot{P},\dot{Q}) := \int_{\mbbT}P\dot{Q}dx - H(P,Q).\] 
	This transforms the system \eqref{eq::systemA} into
	\begin{equation*}\label{eq::stoch-hamiltonian}
		\dot{P} = -\frac{\del H}{\del Q} - \frac{\del \hat{H}}{\del Q}\circ \Phi\xi,\quad \dot{Q} = \frac{\del H}{\del P} + \frac{\del \hat{H}}{\del P}\circ \Phi\xi,
	\end{equation*}
	which is the Hamiltonian system \eqref{eq::HamiltonainIntro}. 
	This Hamiltonian induces the following symplectic structure:
	\begin{equation}\label{eq::induced-symplectic-structure}
		\omega(t) = \int_\mbbT dP\wedge dQ dx.
	\end{equation} 
	The following theorem is due to Hong and Chen in \cite{CheHo2016}.
	\begin{theorem}
		The phase flow of the stochastic Schr\"odinger equation
		\begin{equation*}
			i\pa_t u + \lap u + \abs{u}^2u = u\circ\Phi\xi,
		\end{equation*}
		with $ \xi $ interpreted as in \eqref{eq::snlse}, preserves the symplectic structure
		\begin{equation*}
			\omega(t) = \int_\mbbT d\bar{u}\wedge du dx.
		\end{equation*}
	\end{theorem}
	In the context of numerical schemes, we may write the symplectic form in terms of the discretisation of the coordinates by
	\begin{equation}\label{eq:symplectic-form-disc}
		\omega(u,\bar{u})=\int du_{n+1}\wedge d\bar{u}_{n+1}dx = \int du_n\wedge d\bar{u}_n dx,\quad\forall n\geq0.
	\end{equation}
	\begin{remark}
		The exterior, or wedge, product $ \wedge $ is a product on the exterior algebra of differential $ k- $forms. Importantly, if $ a $ and $ b $ are differential $ k- $forms of degree $ p $ and $ q $, then
		\begin{equation}\label{eq:wedge-propery}
			a\wedge b = (-1)^{pq}b\wedge a
		\end{equation}
	\end{remark}
	
	\section{Symplectic Resonance-Based Runge-Kutta Scheme for Dispersive PDE}\label{sec::2}
	
	In this Section, we present the derivation of resonance-based Runge-Kutta methods in the deterministic case, following the treatment laid out in \cite{M.S2023}. We first rewrite the Duhamel formulation \eqref{eq::intro_SNLS_duhamel} in Fourier space with only the deterministic part
	\begin{equation}\label{eq:duhamel_Fourier}
		\begin{aligned}
			u_k(t)  & = e^{-it k^2 }v_k - i \lam \sum_{k = -k_1 + k_2 + k_3} e^{-it k^2} \\ & \int_0^t e^{is k^2} \bar{u}_{k_1}(s) u_{k_2}(s) u_{k_3}(s)ds ,
		\end{aligned}
	\end{equation}
	where the operator $ e^{it \Delta} $ is assigned in the Fourier space to $ e^{-itk^2} $, and the point-wise product is mapped to the convolution product. Then, one observes the following decomposition.
	\begin{equation} \label{first_approximation}
		u_k(t) =  e^{-it k^2 }v_k + R_k(t).
	\end{equation}
	Here, $ R_k(t) $ is of order $ \mathcal{O}(t) $. One substitutes \eqref{first_approximation} inside the integral of the right-hand side of \eqref{eq:duhamel_Fourier}, in order to obtain
	\begin{equation*}
		u_k(t)   = e^{-it k^2 }v_k - i \lam \sum_{k = -k_1 + k_2 + k_3} e^{-it k^2}\int_0^t e^{is (k^2+k_1^2 - k_2^2 - k_3^2)} ds \bar{v}_{k_1} v_{k_2} v_{k_3}  + \mathcal{O}(t^2).
	\end{equation*}
	The oscillatory integral in the expansion above can be rewritten as 
	\begin{equation}\label{eq:duhamel_eq1}
		-i\lam\sum_{k+k_1=k_2+k_3}\int_0^{ t}e^{-2iskk_1+2isk_2k_3}ds.
	\end{equation}
	In \cite{M.S2023}, to ensure the symplecticity of their schemes, the authors introduce the following approximation.
	\begin{equation}\label{discretisation_low_dom}
		\begin{aligned}
			e^{-2iskk_1 + 2isk_2k_3} &= e^{-2iskk_1}+e^{2isk_2k_3} - 1 + \rbr{e^{-2iskk_1}-1}\rbr{e^{2isk_2k_3} - 1} \\
			&\approx e^{-2iskk_1} + e^{2isk_2k_3} - 1 := \mcK(s;k,k_1,k_2,k_3).
		\end{aligned}
	\end{equation}
	The main idea is that, for obtaining a low-regularity scheme, one wants to avoid discretising the operator $ e^{is P(k_1,k_2,k_3)} $ with $ P(k_1,k_2,k_3) =  -2kk_1 + 2k_2k_3 $
	via a straightforward Taylor expansion. This expansion would require more regularity on the initial datum $v$. 
	The original idea in \cite{O.S2018} was to split the operator into a dominant part $ P_{\text{\tiny{dom}}} $ and a lower part $  P_{\text{\tiny{low}}}  $:
	\begin{equation} \label{naive_splitting}
		P = P_{\text{\tiny{dom}}}  + P_{\text{\tiny{low}}}, \quad P_{\text{\tiny{dom}}}  = 2 k_1^2, \quad P_{\text{\tiny{low}}} =  - 2 k_1 (k_2 + k_3)  + 2 k_2 k_3.
	\end{equation}
	Then, one performs a Taylor expansion on the lower part and integrates the dominant part. The result of the integral can be mapped back to physical space as 
	\begin{equation*}
		\int^t_0 e^{2is k_1^2} ds = \frac{e^{2it k_1^2}-1}{2i k_1^2} \rightarrow \varphi_1(-2 i \Delta t), \quad \varphi_1(z) = \frac{e^z-1}{z}.
	\end{equation*}
	Such a naive splitting, \eqref{naive_splitting}, does not give a symplectic scheme. One has to change it into:
	\begin{equation*}
		P = S_{\text{\tiny{dom}}}  + S_{\text{\tiny{low}}}, \quad S_{\text{\tiny{dom}}}  = - 2 k k_1, \quad S_{\text{\tiny{low}}} =    2 k_2 k_3.
	\end{equation*}
	To preserve symmetries, one has to perform the same treatment for the dominant and the lower part, in the sense that they both need to be integrated exactly:
	\begin{equation*}
		\int^t_0 e^{-2is k k_1} ds = \varphi_1(-2 i t k k_1).
	\end{equation*}
	Then, using the fact that 
	\begin{equation*}
		-2 k k_1 = k^2 + k_1^2 - (k_2 + k_3)^2,
	\end{equation*}
	one has
	\begin{equation*}
		\frac{e^{-2it k k_1}}{k k_1} \bar{v}_{k_1} v_{k_2} v_{k_3} \mapsto \partial_x^{-1} e^{-i t \Delta} \left(  \partial_x^{-1} e^{-i t \Delta} \bar{v}  e^{i t \Delta} v^2 \right).
	\end{equation*}
	One performs the same kind of integration on the lower part. This approach can be generalised to higher orders by introducing the operator $ P_d[f] $, defined as follows. We first fix $d \in \mathbb{N}$ distinct interpolation points $ 0 \leq \gamma_1 \leq \cdots \leq \gamma_d \leq 1 $, then $ P_d[f] $ is the unique interpolation polynomial of degree $d-1$ that matches the function values of $f$ at the points $t \gamma_j $:
	\begin{equation} \label{interpolation}
		P_d[f](t \gamma_j) = f(t  \gamma_j ).
	\end{equation}
	A good approximation of the kernel $ \exp(-2iskk_1 + 2isk_2k_3) $ is 
	\begin{align*}
		\mathcal{K}_{2d}(s; k, k_1, k_2, k_3) & := e^{-2iskk_1}
		P_d[e^{2i \cdot k_2k_3}](s) + e^{2isk_2k_3}P_d[e^{-2i \cdot kk_1}](s)
		\\ & \quad -P_d[e^{2i \cdot k_2k_3}](s) P_d[e^{-2i \cdot kk_1}](s).
	\end{align*}
	From \cite[Proposition 3.6]{M.S2023}, one has
	\begin{equation} \label{bound_kernel}
		\mathcal{K}_{2d}(s; k, k_1, k_2, k_3) - e^{-2iskk_1 + 2isk_2k_3} = \mathcal{O}(t^{d+1} (k k_1 k_2 k_3)^{d/2}).
	\end{equation}
	Let us mention that the term $ \frac{1}{k k_1} $ can be mapped back to physical space only in dimension one. For higher dimensions, this is not true, and one must perform further approximations that require more regularity of the initial datum. The choice \eqref{discretisation_low_dom} guarantees the following symmetry, crucial for having a symplectic scheme:
	\begin{equation}\label{eq:symmetry}
		\mcK_{2d}(s;k,k_1,k_2,k_3)=\overline{\mcK_{2d}(s;k_2,k_3,k,k_1)}.
	\end{equation}
	This property ensures the preservation of the symplectic form, \eqref{eq:symplectic-form}. We first rewrite an approximation of the Duhamel formulation;
	\begin{equation*}
		\begin{aligned}
			u_k(t)  & \approx e^{-it k^2 }v_k - i \lam \sum_{k = -k_1 + k_2 + k_3} e^{-it k^2}
			\\ & \int_0^t 	\mcK_{2d}(s;k,k_1,k_2,k_3) \bar{w}_{k_1}(s) w_{k_2}(s) w_{k_3}(s)ds + \mathcal{O}(t^2 k_1 k_2) 
		\end{aligned}
	\end{equation*}
	where $ w_{k_1}(s) = e^{-is k_1^2 } v_{k_1} $. Then, we want to use a polynomial-type interpolation for the term involving the unknown:
	\begin{equation*}
		w_{k_i}(s) \approx \sum_{p=0}^S \frac{s^p}{t^p} \sum_{q=0} a_{p,q} w_{k_i}( c_q t) ,
	\end{equation*}
	for some  $ S \in \mathbb{N}$, $ 0 \leq c_0 \leq c_1 < \cdots < c_S \leq 1$,  $a_{p,q} \in \mathbb{C}$.
	Using this interpolation, one defines the map,
	\begin{align*}
		v &\mapsto \mcF^{2d}_{p}(t,c_{q};v)\\
		&:=-i\lam\sum_{k\in\mbbZ}e^{ixk}\sum_{k+k_1=k_2+k_3}\frac{1}{t^{p+1}}\int_0^{c_{q}t }  \mcK_{2d}(s;k,k_1,k_2,k_3) s^{p} ds \overline{v}_{k_1} v_{k_2} v_{k_3},
	\end{align*}
	where $ c_q \in [0,1] $.
	This map, which describes the discretisation of the integral in \eqref{eq:duhamel_eq1}, leads to the following scheme for the NLSE,
	\begin{equation}\label{eq:scheme_nlse}
		\begin{aligned}
			u_{n+1} &= e^{it\pa_x^2}u_n+t e^{it\pa_x^2}\sum_{\alpha\in D}b^{\alpha}K_{\alpha},\\
			K_{\alpha} & =\mathfrak{F}^{2d}_p\rbr{t;c_{q};u_n+t\sum_{\ti\alpha\in D}a^{\ti\alpha}_{\alpha}K_{\ti\alpha}}, \quad
			D = \set{(p,q,r):0 \leq p,q,r \leq S}.
		\end{aligned}
	\end{equation}
	where $ \alpha =(p,q,r) $ and $ \mathfrak{F}^{2d}_p $ is the physical space map of $ \mcF^{2d}_{p} $.
	Under the correct constraints on the parameters $ a^\alpha_{\ti\alpha}$ and $b^{\alpha} $, this class of schemes is capable of conserving the Hamiltonian and the $ L^2 $ norm. We quote \cite[Theorem 4.1]{M.S2023}.
	\begin{theorem} \label{condition_deterministic} Suppose that  $ u \in H^r $, $r > 3/2$, and that the real-valued coefficients  $ b^{\alpha}$, $ a^{\tilde{\alpha}}_{\alpha} $
		satisfy
		\begin{equation*}
			b^{\ti\alpha}b^{\alpha}=b^{\alpha}a^{\ti\alpha}_{\alpha}+b^{\ti\alpha}a^{\alpha}_{\ti\alpha},\quad\forall \alpha,\ti\alpha\in S,
		\end{equation*}
		Then the Runge–Kutta resonance-based schemes \eqref{eq:scheme_nlse}
		preserve the corresponding quadratic first integrals exactly.
	\end{theorem} 
	
	\section{Stochastic Symplectic Resonance-Based Runge-Kutta Method}\label{sec::4}
	We now explain how to derive the midpoint rule for the SNLSE; the idea is to adapt the general scheme given in \cite{B.S2022}. 
	We recall the Duhamel formulation, given by \eqref{eq::snlse} in physical space, and in Fourier space by,
	\begin{equation}\label{eq:duhamel_Fourier_stochastic}
		\begin{aligned}
			u_k(t)  & = e^{-it k^2 }v_k - i \lam \sum_{k = -k_1 + k_2 + k_3} e^{-it k^2}\int_0^t e^{is k^2} \bar{u}_{k_1}(s) u_{k_2}(s) u_{k_3}(s)ds
			\\ &  \quad - i \kappa \sum_{k = k_1 +k_2} e^{-it k^2}\int_0^t e^{is k^2} u_{k_1}(s) \Phi_{k_2} \circ dW_{k_2}(s). 
		\end{aligned}
	\end{equation}
	Iterating the Duhamel formula and truncating at order 
	$ \mathcal{O}(t^{\frac{3}{2}}) $, one gets
	\begin{align*}
		u_k(t)  & = e^{-it k^2 }v_k - i \lam \sum_{k = -k_1 + k_2 + k_3} e^{-it k^2}\int_0^t e^{is (k^2+k_1^2 - k_2^2 - k_3^2)} ds \bar{v}_{k_1} v_{k_2} v_{k_3} 
		\\ &	\quad - i \kappa
		\sum_{k = k_1 +k_2} e^{-it k^2}\int_0^t e^{is k^2} e^{-i s k_1^2}v_{k_1}\Phi_{k_2} \circ  dW_{k_2}(s)
		\\ & \quad- \kappa^2 \sum_{k = k_1 +k_2 +k_3} e^{-it k^2}\int_0^t e^{is k^2} \\&  e^{-is (k_1 + k_2)^2} \int_0^s e^{ir (k_1+k_2)^2} v_{k_1}\Phi_{k_2}\Phi_{k_3} \circ dW_{k_2}(r) \circ  dW_{k_3}(s)
		+ \mathcal{O}(t^{\frac{3}{2}}).
	\end{align*}
	In the next subsection, we explain how to discretise the stochastic integral of order $ \sqrt{t} $ and show that the integral of order $ t $ disappears up to order $ t^2 $.
	\subsection{Discretisation of Stochastic Iterated Integrals}
	\label{discretisation_stochastic_integral}
	The first stochastic integral is given in Fourier space by
	
	\begin{equation}\label{eq:stoch_integral_1}
		I_1 = -i \kappa  \sum_{k\in\mbbZ}e^{ikx}\sum_{k =k_1+k_2}\Phi_{k_2}v_{k_1}e^{-it(k_1+k_2)^2}\int_0^{t}e^{is(kk_2+k_1k_2)}\circ dW_{k_2}(s).
	\end{equation}
	One way of approaching the discretisation of this integral is to perform integration by parts to rewrite the right-hand side as a random boundary term and a Riemann integral over a Brownian path. The Riemann integral can then be approximated via a quadrature rule. For instance, we could write
	\[
	\int_{0}^{t} e^{i s P(k)} dW(s)
	\;=\;
	e^{ itP(k)}\,W(t)
	-
	i P(k)\int_{0}^{t} e^{i s P(k)}\,W(s)\,\mathrm ds,
	\]
	and approximate this by 
	\[
	\int_{0}^{t} e^{is P(k)}\,W(s) ds
	=
	\frac{t}{2}\Bigl( W(0) + e^{ i t P(k)}W(t) \Bigr)
	+\mathcal O\!\bigl(t^{3/2}\bigr).
	\]
	It is not immediately obvious whether it offers any advantage or whether this discretisation will preserve relevant structures. For now, we proceed with a simple Taylor expansion of the operator, which suffices to obtain an improvement in the regularity of initial data,
	\begin{equation}\label{eq:stoch_integral_taylor}
		\begin{aligned}
			I_1 &= -i\kappa \sum_{k \in \mathbb{Z}} e^{ikx} \sum_{\substack{k_1, k_2 \in \mathbb{Z} \\ k = k_1 + k_2}} \Phi_{k_2} v_{k_1} e^{-it(k_1 + k_2)^2} \int_0^t e^{is(kk_2 + k_1k_2)} \circ dW_{k_2}(s) \\
			&\approx -i\kappa \sum_{k \in \mathbb{Z}} e^{ikx} \sum_{\substack{k_1, k_2 \in \mathbb{Z} \\ k = k_1 + k_2}} \Phi_{k_2} v_{k_1} e^{-it(k_1 + k_2)^2} \int_0^t (1 + is(kk_2 + k_1k_2) \\ 
			& \hspace{180pt}+ \mathcal{O}(s^2(kk_2 + k_1k_2)^2))\circ dW_{k_2}(s) \\
			&= -i\kappa \sum_{k \in \mathbb{Z}} e^{ikx} \sum_{\substack{k_1, k_2 \in \mathbb{Z} \\ k = k_1 + k_2}} \Phi_{k_2} v_{k_1} e^{-it(k_1 + k_2)^2} \int_0^t \circ dW_{k_2}(s) \\
			&\quad - \kappa \sum_{k \in \mathbb{Z}} e^{ikx} \sum_{\substack{k_1, k_2 \in \mathbb{Z} \\ k = k_1 + k_2}} \Phi_{k_2} v_{k_1} e^{-it(k_1 + k_2)^2} \int_0^t is(kk_2 + k_1k_2) \circ dW_{k_2}(s) \\
			&\quad + \mathcal{O}\left(\kappa \sum_{k \in \mathbb{Z}} \sum_{\substack{k_1, k_2 \in \mathbb{Z} \\ k = k_1 + k_2}} k_1k_2^2 t^{3/2}\right) \quad a.s.
		\end{aligned}
	\end{equation}
	For a structure-preserving approximation, one has to replace $ u_{k_1} $ by a midpoint approximation:
	\begin{equation*}
		v_{k_1}:=	u_{k_1}(0) \rightarrow \frac{u_{k_1}(0) + e^{i t k_1^2}u_{k_1}(t) }{2}, 
	\end{equation*}
	or a similar symmetric approximation at various points. In the end, one gets the following almost sure approximation:
	\begin{equation*}
		I_1^{S} = -i \kappa \sum_{k\in\mbbZ}e^{ikx}\sum_{k=k_1+k_2}\Phi_{k_2} \frac{u_{k_1}(0) + e^{i t k_1^2}u_{k_1}(t) }{2} e^{-it(k_1+k_2)^2} \int_0^t \circ d W_{k_2}(s).
	\end{equation*}
	In a general setting, one can consider the map $ \mcP $ defined by
	\begin{equation}\label{eq:general-disc-map-stoch}
		\begin{aligned}
			\mcP_p^{2d}(t;c_q;v) & = -i  \kappa\sum_{k\in\mbbZ}e^{ixk}\frac{1}{t^{p+\frac12}}\sum_{k = k_1 + k_2}v_{k_1}\Phi_{k_2} \\ & \int^{c_q t}_0 P_d[e^{i \cdot (kk_2 + k_1k_2)}](s) s^{p}  \circ dW_{k_2}(s).
		\end{aligned}
	\end{equation}
	For the double-iterated stochastic integral $I_2$, one has
	\begin{align*}
		I_2 & = 	- \kappa^2 \sum_{k = k_1 +k_2 +k_3} e^{-it k^2}\int_0^t e^{is k^2}  \\ & e^{-i s(k_1+k_2)^2} v_{k_1}\Phi_{k_2}\Phi_{k_3}\int_0^s e^{ir (k_1+k_2)^2} e^{-ir k_1^2} \circ dW_{k_2}(r) \circ dW_{k_3}(s).
	\end{align*}
	We proceed by expansion of both dominant and lower terms to get the following discretisation:
	\begin{equation*}
		I_2 \approx	-\kappa^2 \sum_{k = k_1 +k_2 +k_3}v_{k_1} \Phi_{k_2}\Phi_{k_3} \int_0^t   \int_0^s  \circ  dW_{k_2}(r) \circ  dW_{k_3}(s)+ \mathcal{O}(k_1k^2_2 k_3^2t^{\frac{3}{2}}).
	\end{equation*}
	This discretisation was proposed in \cite{AB2023}, but it will not preserve the structure. One has to make use of the key idea in \cite{B.B.M.S2023}, which is to iterate the Duhamel formula symmetrically, for example, by using a midpoint rule. Then, $I_2$ is transformed into
	\begin{align*}
		I_2 &=	- \kappa^2 \sum_{k = k_1 +k_2 +k_3} e^{-it k^2}\int_0^t e^{is (k^2- (k_1 + k_2)^2)}   \left( \frac{1}{2} v_{k_1}\Phi_{k_2} \int_0^s e^{ir (k_1+k_2)^2}   \right.
		\\ & \left. e^{-irk_1^2}  \circ dW_{k_2}(r)	+ \frac{1}{2}    v_{k_1} \Phi_{k_2}\Phi_{k_3} \int_t^s e^{ir (k_1+k_2)^2} e^{-irk_1^2} \circ  dW_{k_2}(r) \right)\circ dW_{k_3}(s).
	\end{align*}
	Then, proceeding with a full Taylor expansion on the various operators, one gets
	\begin{align*}
		I_2 &= -\kappa^2\sum_{k\in\mbbZ}e^{ikx}\sum_{k = k_1 +k_2 +k_3}e^{itk} \\
		&\quad v_{k_1}\Phi_{k_2}\Phi_{k_3}\int_0^t  \bigg( \frac{1}{2}\int_0^s\circ dW_{k_2}(r)  + \frac12\int_t^s \circ dW_{k_2}(r)\bigg) \circ dW_{k_3}(s) \\
		&\quad+\mcO(k_1k_2^2k_3^2t^{2})\quad a.s.
	\end{align*}
	We observe that 
	\begin{equation*}
		\frac{1}{2}\int_0^s\circ dW_{k_2}(r)  + \frac12\int_t^s \circ dW_{k_2}(r) = W_{k_2}(s) - \frac12W_{k_2}(t),
	\end{equation*}
	and so, integrating against $ W_{k_3}(t) $, we obtain
	\begin{equation*}
		\int_0^t \rbr{W_{k_2}(s) - \frac12W_{k_2}(t)} \circ dW_{k_3}(s) = \int_0^t W_{k_2}(s) \circ dW_{k_3}(s) - \frac12W_{k_2}(t)W_{k_3}(t).
	\end{equation*}
	But in Stratonovich theory, we have the following,
	\begin{equation*}
		\int_0^t W_{k_2}(s) \circ dW_{k_3}(s) + \int_0^t W_{k_3}(s) \circ dW_{k_2}(s) = W_{k_2}(t)W_{k_3}(t),
	\end{equation*} which implies that
	\begin{equation}\label{eq:order-of-iterate}
		I_2 = 0 + \mcO(k_1k_2^2k_3^2t^2) \quad a.s.
	\end{equation}
	In short, when constructing first-order schemes for \eqref{eq::snlse} with Stratonovich-type noise, it makes sense to disregard the second iteration and to set the approximation of $ I_2 $ to be zero. This simplification has been observed in \cite{B.B.M.S2023} for double-iterated integrals in the deterministic setting.
	\subsection{Discretisation Maps Properties}
	When verifying that our scheme preserves certain structures, we will need several hypotheses. These are based on formal calculations that are performed when deriving the relevant properties.
	\begin{hypothesis}\label{hyp:Phi}
		The smoothing, or covariance, operator of the noise is a self-adjoint Hilbert-Schmidt operator $ \Phi:L^2\rightarrow H^2 $, which is symmetric in its Fourier coefficients, i.e.
		\begin{equation*}
			\Phi^*_k = \Phi_k,\quad\Phi_k = \Phi_{-k}.
		\end{equation*}	
	\end{hypothesis}
	\begin{remark}
		An example of $\ \Phi $ with the required properties is given by,
		\begin{equation*}
			\Phi_k =  \begin{cases}
				0,\quad k=0, \\
				\frac{1}{k^2},\quad k\in\mbbZ\backslash \set{0},	
			\end{cases}
		\end{equation*}
		and this indeed provides the required regularity to the noise if we consider first-order Taylor expansions of the kernel. 
	\end{remark}
	
	\begin{scheme}\label{scheme:stoch-runge-kutta}
		We introduce the following scheme as the stochastic extension of the Runge-Kutta resonance-based scheme for the NLSE:
		\begin{align*}
			u_{n+1} &= e^{it \pa_x^2} u_n + t e^{it \pa_x^2} \sum_{\alpha \in D} b^{(0)}_{\alpha} K_{\alpha} + \sqrt{t} e^{it \pa_x^2} \sum_{\alpha \in D} b^{(1)}_{\alpha} L_\alpha, \\
			U_{\alpha}(u_n) &= u_n + t\sum_{\ti{\alpha} \in D} a^{(0)}_{\alpha,\ti{\alpha}} K_{\ti{\alpha}} + \sqrt{t}\sum_{\ti{\alpha}\in D} a^{(1)}_{\alpha,\ti{\alpha}} L_{\ti{\alpha}},
		\end{align*}
		where 	$D = \set{(p,q,r):0 \leq p,q,r \leq S},$ and 
		\begin{equation*}
			K_{\alpha} = \mathfrak{F}^{2d}_{p}\rbr{t;c_{q};U_{\alpha}(u_n)},\quad
			L_{\alpha} = \mathfrak{P}^{2d}_{p}\rbr{t;c_{q};U_{\alpha}(u_n)}.\
		\end{equation*}
		Here $ \alpha =(p,q,r) $ and $ \mathfrak{P}^{2d}_p $ is the map $ \mathcal{P}^{2d}_p $ in physical space.
	\end{scheme}
	We will show that, similarly to in the deterministic setting, the conditions for norm preservation will be an algebraic system of polynomials in terms of the coefficients $ a^{(i)}_{\alpha, \ti{\alpha}}$ and $b^{(i)}_{\alpha}$ for $ \alpha\in D $ which read as follows:
	\begin{hypothesis}\label{hyp:coefficient}
		The coefficients of Scheme \ref{scheme:stoch-runge-kutta}, are such that
		\begin{equation*}
			b_\alpha^{(i)}b_{\ti{\alpha}}^{(j)} - b^{(i)}_{\alpha} a^{(j)}_{\alpha,\ti{\alpha}} - b^{(j)}_{\ti{\alpha}} a^{(i)}_{\ti\alpha,\alpha}=0,\quad i,j=0,1.
		\end{equation*}
	\end{hypothesis}
	We now introduce an explicit example of the above class of schemes, the stochastic resonance-based midpoint rule.
	\begin{example}[Stochastic Resonance-Based Midpoint Rule]\label{ex:stoch-res-midpoint}
		We consider the case $ S=0 $, with the choice $c_0 = 0$, $ b^{(i)}_{\alpha} = 1 $ and $ a^{(i)}_{\alpha, \alpha}  = \frac12 $ with $\alpha = (0,0,0)$, and then we get the scheme
		\begin{align*}
			u_{n+1} &= e^{it\pa_x^2}u_n + t  e^{it\pa_x^2}K_{0,0,0} + \sqrt{t} e^{it\pa_x^2}L_{0,0,0},
			\\
			U_{0,0,0}(u_n) &= u_n +  \frac{t}{2} K_{0,0,0} + \frac{\sqrt{t}}{2}  L_{0,0,0},
			\\
			K_{0,0,0} &= \mathfrak{F}_{0}\rbr{t;1;U_{0,0,0}(u_n)},\quad
			L_{0,0,0} = \mathfrak{P}_{0}\rbr{t;1;U_{0,0,0}(u_n)}.\
		\end{align*}
		From the first identity, one observes that
		\begin{equation*}
			e^{-it\pa_x^2}u_{n+1}-u_n =  t  K_{0,0,0} + \sqrt{t} L_{0,0,0}.
		\end{equation*}
		Therefore, 
		\begin{align*}
			U_{0,0,0}(u_n) & = u_n +  \frac{t}{2} K_{0,0,0} + \frac{\sqrt{t}}{2}  L_{0,0,0}
			\\ & = u_n + \frac{e^{-it\pa_x^2}u_{n+1}-u_n}{2} = \frac{u_n + e^{-it\pa_x^2}u_{n+1}}{2}.
		\end{align*}
		We obtain the simplification
		\begin{equation*}
			u_{n+1} = e^{it\pa_x^2}u_n + e^{it\pa_x^2}
			\left(   t \mathfrak{F}_{0}\rbr{t;1;u_{n+\frac12}} +
			\sqrt{t} \mathfrak{P}_{0}\rbr{t;1;u_{n+\frac12}} \right),
		\end{equation*}
		where
		\begin{equation*}
			u_{n+\frac12} = \frac{u_n + e^{-it\pa_x^2}u_{n+1}}{2}.
		\end{equation*}
		This is the stochastic resonance-based midpoint rule. By expanding the discretisation maps $ \mathfrak{F}_{0}, \mathfrak{P}_{0}$ we obtain
		\begin{align}\label{eqn:resonance_based_midpoint_rule_u_NLSE}
			&	u_{n+1}=e^{i\partial_x^2 t}u_{n}-i\lambda \frac{\lambda}{2}\partial_x^{-1}\left(\left[e^{-i t \partial_x^2}\overline{\partial_x^{-1}u_{n+\frac{1}{2}}}\right]\left[e^{i t \partial_x^2}\left(u_{n+\frac{1}{2}}\right)^2\right]\right)\notag\\ & -\frac{\lambda}{2}e^{i\partial_x^2 t}\partial_x^{-1}\left(\overline{\partial_x^{-1}u_{n+\frac{1}{2}}}\left(u_{n+\frac{1}{2}}\right)^2\right) + \frac{\lambda}{2} e^{i t\partial_x^2 }\overline{u_{n+\frac{1}{2}}}e^{-i\tau \partial_x^2}\left(e^{i t \partial_x^2}\partial_x^{-1}u_{n+\frac{1}{2}}\right)^2\notag \\ & -\frac{\lambda}{2} e^{i t\partial_x^2 } \overline{u_{n+\frac{1}{2}}}\left(\partial_x^{-1}u_{n+\frac{1}{2}}\right)^2 + i\lambda e^{i t \partial_x^2 }  \tau|u_{n+\frac{1}{2}}|^2u_{n+\frac{1}{2}}\notag \\
			&\quad-i\lambda t \left[\int_{\mathbb{T}}|u_{n+\frac{1}{2}}|^2u_{n+\frac{1}{2}}dx +\overline{u^{n+\frac{1}{2},0}}e^{i t\partial_x^2 }(u_{n+\frac{1}{2}})^2-\overline{u_{n+\frac{1}{2},0}}\int_{\mathbb{T}}\left(u_{n+\frac{1}{2}}\right)^2 dx \right. \notag \\ & \left.  +  2 u_{n+\frac{1}{2},0}e^{i t\partial_x^2 }\left(|u_{n+\frac{1}{2}}|^2\right) -\left(u_{n+\frac{1}{2},0}\right)^2e^{it\partial_x^2 }\overline{u_{n+\frac{1}{2}}}\right]
			\notag \\ &  -i \kappa  e^{it \partial_x^2} \left(  u_{n+\frac{1}{2}} \int_0^t \circ \Phi d W(s) \right),
		\end{align}
		where $ u_{n+\frac{1}{2},0} $ is the zero Fourier mode of $  u_{n+\frac{1}{2}}   $.
		The deterministic part of the scheme was obtained in \cite[Example 4.6]{M.S2023}
	\end{example}
	
	\section{Conserved Quantities}\label{sec::5}
	In this section, we show that in the presence of conservative noise, assuming Lemma \ref{prop:disc-maps} holds, the stochastic resonance Runge-Kutta Scheme \ref{scheme:stoch-runge-kutta} preserves both the symplectic structure and 
	the mass of the SNLSE \eqref{eq::snlse} under Hypotheses \ref{hyp:Phi}, \ref{hyp:coefficient}. This is true when the noise is multiplicative, has real Fourier coefficients, and is symmetric in Fourier space.   
	\begin{theorem}\label{thm:Symplecticity}
		Under Hypotheses \ref{hyp:Phi}, \ref{hyp:coefficient}, and assuming Lemma \ref{prop:disc-maps} holds, the stochastic resonance Runge-Kutta Scheme \ref{scheme:stoch-runge-kutta} is symplectic.
	\end{theorem}
	Before proving Theorem \ref{thm:Symplecticity} by breaking it into two propositions, we need the following lemma:
	\begin{lemma}\label{prop:disc-maps}
		The discretisation maps $ \mathfrak{F}, \mathfrak{P} $ of the numerical Scheme \ref{scheme:stoch-runge-kutta} are such that
		\begin{equation}\label{eq:orthogonality}
			\Re\sqbr{\int_{\mbbT}\overline{u} \mathfrak{F}^{2d}_{p}\rbr{t;c_{q};u}  dx} = 0, \quad \Re\sqbr{\int_{\mbbT}\overline{u} \mathfrak{P}^{2d}_{p}\rbr{t;c_{q};u}  dx} = 0.
		\end{equation}	
		where $ \Re(f) $ represents the real part of $ f $.
	\end{lemma}
	\begin{proof}
		For $\mathfrak{F}$, one has 
		\begin{equation*}
			\Re\sqbr{\sum_{k\in\mbbZ}\rbr{-i\lambda}\sum_{k+k_1=k_2+k_3} \mcK_{2d}(s;k,k_1,k_2,k_3)\overline{u_k u_{k_1}}u_{k_2}u_{k_3}}=0.
		\end{equation*}
		This is due to the symmetry
		\begin{equation*}
			\mcK_{2d}(s;k,k_1,k_2,k_3)=\overline{\mcK_{2d}(s;k_2,k_3,k,k_1)}.
		\end{equation*}
		Multiplying by $ s^pt^{-p+1} $ and integrating over $ [0,c_qt] $ results in the first equation in \eqref{eq:orthogonality}.
		For $ \mathfrak{P}^{2d}_p $, one has
		\begin{align*}
			A &= \int_\mbbT \bar{u}\mathfrak{P}^{2d}_p(t;c_q;u)dx \\
			& = -i  \kappa\frac{1}{t^{p+1}}\sum_{k = k_1 + k_2} \bar{u}_k u_{k_1}  \int^{c_q t}_0 P_d[e^{i \cdot (kk_2 + k_1k_2)}](s) s^{p} \Phi_{k_2}\circ dW_{k_2}(s).
		\end{align*}
		We note that $ W_{k_2}(s) $ takes values in $\mathbb{R}$ and $ \Phi_{k_2} = \bar{\Phi}_{k_2} $, by Hypothesis \ref{hyp:Phi}, which implies that
		\begin{equation*}
			\overline{A} = 	 i  \kappa\frac{1}{t^{p+1}}\sum_{ k_1 = k- k_2} \bar{u}_k u_{k_1}  \int^{c_q t}_0 P_d[e^{-i \cdot (kk_2 + k_1k_2)}](s) s^{p} \Phi_{k_2}\circ dW_{k_2}(s).
		\end{equation*} 
		Then, one can perform the change of variable $ k_2' = -k_2 $ to see that $ \overline{A} = -A $ which allows us to conclude.
	\end{proof}
	\begin{proposition}\label{prop:conservation-of-mass}
		Under Hypothesis \ref{hyp:coefficient} and Lemma \ref{prop:disc-maps}, the Runge-Kutta stochastic resonance scheme \ref{scheme:stoch-runge-kutta} conserves the mass of \eqref{eq::snlse}, that is,
		\begin{equation*}
			\int_\mbbT \abs{u_{n+1}}^2 dx = \int_\mbbT \abs{u_{n}}^2 dx.
		\end{equation*}		
	\end{proposition}
	\begin{proposition}\label{prop:conservation-of-energy}
		Under Hypotheses \ref{hyp:Phi} and \ref{hyp:coefficient}, the stochastic resonance Runge-Kutta Scheme \ref{scheme:stoch-runge-kutta} conserves the symplectic structure \eqref{eq::induced-symplectic-structure} induced by the generalised Hamiltonian,
		\begin{align*}
			\dot{u} = -\frac{\delta H}{\delta \bar{u}} - \frac{\delta \hat{H}}{\delta \bar{u}} \circ  \Phi\xi, \quad
			\dot{\bar{u}} = \frac{\delta H}{\delta u} + \frac{\delta \hat{H}}{\delta u} \circ  \Phi\xi,
		\end{align*}
		where
		\begin{equation*}
			H(u) =  \frac12\int\abs{\grad u}^2 dx + \frac{\lam}{4}\int \abs{u}^{4} dx, \quad
			\hat{H}(u) = \frac{\kappa}{2}\int \abs{u}^2 dx,
		\end{equation*}
		of the equation \eqref{eq::snlse}.
	\end{proposition}
	\begin{proof}[Proof of Proposition \ref{prop:conservation-of-mass}]
		We wish to show that $ \int \abs{u_{n+1}}^2dx = \int \abs{u_{n}}^2dx $. First we note that, since $ e^{it\pa_x^2} $ is unitary we have
		\begin{align*}\label{eq:mass-conservation-id}
			\abs{u_{n+1}}^2 - \abs{u_{n}}^2 &= \abs{u_{n+1}}^2 - \abs{e^{it\pa_x^2}u_{n}}^2 \\
			&=\abs{u_{n+1}-e^{it\pa_x^2}u_n}^2 + 2\Re\rbr{\overline{e^{it\pa_x^2}u_n}\rbr{u_{n+1}-e^{it\pa_x^2}u_n}}\\
			&=: \abs{A}^2 + 2\Re{B},
		\end{align*}
		and introduce the shorthand notation:
		\begin{equation*}
			K'_\alpha = e^{it\pa_x^2} K_\alpha, \quad L'_\alpha = e^{it\pa_x^2} L_\alpha.
		\end{equation*} 
		Recall that from the definition of the Scheme \ref{scheme:stoch-runge-kutta},
		\begin{equation}
			u_{n+1} -e^{it \pa_x^2} u_n =  t e^{it \pa_x^2} \sum_{\alpha \in D} b^{(0)}_{\alpha} K_{\alpha} + \sqrt{t} e^{it \pa_x^2} \sum_{\alpha \in D} b^{(1)}_{\alpha} L_\alpha. 
		\end{equation}
		Hence, we have:
		\begin{align*}
			\abs{A}^2 &=\abs{u_{n+1}-e^{it\pa_x^2}u_n}^2 \\&=\bigg|t\sum_{\alpha \in D} b^{(0)}_{\alpha}e^{it\pa_x^2}K_{\alpha} + \sqrt{t}\sum_{\alpha \in D} b^{(1)}_{\alpha}e^{it\pa_x^2}L_\alpha\bigg|^2 \\
			&=\sum_{\alpha,\ti{\alpha} \in D}t^2 K'_\alpha\overline{K'_{\ti{\alpha}}}b_\alpha^{(0)}b_{\ti{\alpha}}^{(0)} + tL'_\alpha\overline{L'_{\ti{\alpha}}}b_\alpha^{(1)}b_{\ti{\alpha}}^{(1)} 
			+t^{3/2} K'_\alpha\overline{L'_{\ti{\alpha}}}b_\alpha^{(0)}b_{\ti{\alpha}}^{(1)} + t^{3/2}L'_\alpha\overline{K'_{\ti{\alpha}}}b_\alpha^{(1)}b_{\ti{\alpha}}^{(0)}.
		\end{align*}
		Now, following from Lemma \ref{prop:disc-maps} we observe that,
		\begin{equation}\label{eq::conj-rel}
			\begin{aligned}
				\mfR\sqbr{\int_{\mbbT}(\overline{v_1}+\overline{v_2}+\overline{v_3})\mathfrak{F}^{2d}_{p}\rbr{t;c(q);v_1 + v_2 + v_3}}&=0, \\
				\mfR\sqbr{\int_{\mbbT}(\overline{v_1}+\overline{v_2}+\overline{v_3})\mathfrak{P}^{2d}_{p}\rbr{t;c(q);v_1 + v_2 + v_3}}&=0 .
			\end{aligned} 
		\end{equation}
		Set $ v_1 = e^{it\pa_x^2}u_n $,
		$ v_2 = t\sum_{\ti{\alpha}\in D} a^{(0)}_{\alpha,\ti{\alpha}}K'_{\ti{\alpha}}$ and $v_3= \sqrt{t}\sum_{\ti{\alpha}\in D} a^{(1)}_{\alpha,\ti{\alpha}}L'_{\ti{\alpha}}$. Then we have
		\begin{align*}
			\int_{\mbbT}\Re(B)dx &=\int_{\mbbT}\Re\rbr{\overline{e^{it\pa_x^2}u_n}\rbr{u_{n+1}-e^{it\pa_x^2}u_n}}dx \\&=\int_{\mbbT}\Re\rbr{\overline{v_1}\rbr{u_{n+1}-e^{it\pa_x^2}u_n}}dx\\&= \int_{\mbbT}\Re\bigg(\overline{v_1}\bigg(  t\sum_{\alpha \in D} b^{(0)}_{\alpha} K'_{\alpha} + \sqrt{t}\sum_{\alpha \in D} b^{(1)}_{\alpha} L'_\alpha\bigg)\bigg) dx\\
			&=t\sum_{\alpha \in D}b^{(0)}_{\alpha}\Re\int_{\mbbT}\overline{v_1}\mathfrak{F}^{2d}_{p}\rbr{t;c(q);v_1 + v_2 + v_3}dx \\ &\quad +  \sqrt{t} \sum_{\alpha \in D} b^{(1)}_{\alpha}\Re\int_\mbbT\overline{v_1} \mathfrak{P}^{2d}_{p}\rbr{t;c(q);v_1 + v_2 + v_3} dx\\
			&=-t\sum_{\alpha \in D}b^{(0)}_{\alpha}\Re\int_{\mbbT}\overline{v_2+v_3}\mathfrak{F}^{2d}_{p}\rbr{t;c(q);v_1 + v_2 + v_3}dx \\ &\quad -  \sqrt{t} \sum_{\alpha \in D} b^{(1)}_{\alpha}\Re\int_\mbbT\overline{v_2+v_3} \mathfrak{P}^{2d}_{p}\rbr{t;c(q);v_1 + v_2 + v_3} dx \\
			&= -t\sum_{\alpha,\ti\alpha\in D}b_\alpha^{(0)}\Re\int_\mbbT \rbr{ta^{(0)}_{\alpha,\ti\alpha}\overline{K'_{\ti\alpha}}+\sqrt{t}a^{(1)}_{\alpha,\ti\alpha}\overline{L'_{\ti\alpha}}}K'_\alpha dx\\
			&\quad -\sqrt{t}\sum_{\alpha,\ti\alpha\in D}b_\alpha^{(1)}\Re\int_\mbbT \rbr{ta^{(0)}_{\alpha,\ti\alpha}\overline{K'_{\ti\alpha}}+\sqrt{t}a^{(1)}_{\alpha,\ti\alpha}\overline{L'_{\ti\alpha}}}L'_\alpha dx.
		\end{align*}
		This implies that we have
		\begin{align*}
			2\int_{\mathbb{T}} \Re(B) \, dx 
			= -\sum_{\alpha, \tilde{\alpha} \in D} \bigg[ 
			&t^2 \big( b^{(0)}_\alpha a^{(0)}_{\alpha, \tilde{\alpha}} 
			+ b^{(0)}_{\tilde{\alpha}} a^{(0)}_{\tilde{\alpha}, \alpha} \big)
			\Re \int_{\mathbb{T}} \overline{K'_{\tilde{\alpha}}} K'_\alpha \, dx \\
			&+ t^{3/2} \big( b^{(0)}_\alpha a^{(1)}_{\alpha, \tilde{\alpha}} 
			+ b^{(1)}_{\tilde{\alpha}} a^{(0)}_{\tilde{\alpha}, \alpha} \big)
			\Re \int_{\mathbb{T}} \overline{L'_{\tilde{\alpha}}} K'_\alpha \, dx \\
			&+ t^{3/2} \big( b^{(1)}_\alpha a^{(0)}_{\alpha, \tilde{\alpha}} 
			+ b^{(0)}_{\tilde{\alpha}} a^{(1)}_{\tilde{\alpha}, \alpha} \big)
			\Re \int_{\mathbb{T}} \overline{K'_{\tilde{\alpha}}} L'_\alpha \, dx \\
			&+ t \big( b^{(1)}_\alpha a^{(1)}_{\alpha, \tilde{\alpha}} 
			+ b^{(1)}_{\tilde{\alpha}} a^{(1)}_{\tilde{\alpha}, \alpha} \big)
			\Re \int_{\mathbb{T}} \overline{L'_{\tilde{\alpha}}} L'_\alpha \, dx 
			\bigg].
		\end{align*}
		Hence, by combining the two calculations, we get
		\begin{align*}
			\int_{\mbbT}\abs{u_{n+1}}^2 - \abs{u_n}^2 dx &= \int_{\mbbT}\abs{u_{n+1} - e^{it\pa_x^2}u_n}^2 - 2\Re\rbr{\overline{e^{it\pa_x^2}u_n}(u_{n+1} - e^{it\pa_x^2}u_n)}dx \\
			&=\sum_{\alpha,\ti{\alpha} \in D}t^2\Re\int_\mbbT K'_\alpha\overline{K'_{\ti{\alpha}}} \, dx\rbr{b_\alpha^{(0)}b_{\ti{\alpha}}^{(0)} - b^{(0)}_{\alpha} a^{(0)}_{\alpha,\ti{\alpha}} - b^{(0)}_{\ti{\alpha}} a^{(0)}_{\ti\alpha,\alpha}} \\
			&\quad+ t\Re\int_\mbbT L'_\alpha\overline{L'_{\ti{\alpha}}} \, dx\rbr{b_\alpha^{(1)}b_{\ti{\alpha}}^{(1)} - b^{(1)}_{\alpha}a^{(1)}_{\alpha,\ti{\alpha}} - b^{(1)}_{\ti\alpha}a^{(1)}_{\ti\alpha,{\alpha}}} \\
			&\quad+ t^{3/2} \Re \int_\mbbT  K'_\alpha\overline{L'_{\ti{\alpha}}} \, dx \rbr{b_\alpha^{(0)}b_{\ti{\alpha}}^{(1)} - b^{(0)}_{\alpha}a^{(1)}_{\alpha,\ti{\alpha}} - b^{(1)}_{\ti\alpha}a^{(0)}_{\ti\alpha,{\alpha}} }\\
			&\quad+ t^{3/2}\Re\int_\mbbT  L'_\alpha\overline{K'_{\ti{\alpha}}} \, dx \rbr{b_\alpha^{(1)}b_{\ti{\alpha}}^{(0)} - b^{(1)}_{\alpha}a^{(0)}_{\alpha,\ti{\alpha}} - b^{(0)}_{\ti\alpha}a^{(1)}_{\ti\alpha,{\alpha}}}.
		\end{align*}
	\end{proof}
	Our goal is now to show that this scheme is symplectic under the Hypotheses \eqref{hyp:Phi} and \eqref{hyp:coefficient}. To do this, we first want to write the scheme in terms of the Hamiltonian system whose particles are $ \bmu = \set{v_k}_{k\in\mbbZ} $, $ \bnu = \set{\overline{v_k}}_{k\in\mbbZ} $. To do this, we first denote
	\begin{align*}
		\boldsymbol{U}_{\alpha}^{(1)}(\bmu) & = \bmu + t\sum_{\alpha \in D} a^{(0)}_{\alpha,\ti{\alpha}}\bK^{(1)}_{\ti{\alpha}} + \sqrt{t}\sum_{\alpha \in D} a^{(1)}_{\alpha,\ti{\alpha}}\bL^{(1)}_{\ti{\alpha}},
		\\ \boldsymbol{U}_{\alpha}^{(2)}(\bnu) & = \bnu + t\sum_{\alpha \in D} a^{(0)}_{\alpha,\ti{\alpha}}\bK^{(2)}_{\ti{\alpha}} + \sqrt{t}\sum_{\alpha \in D} a^{(1)}_{\alpha,\ti{\alpha}}\bL^{(2)}_{\ti{\alpha}}.
	\end{align*} 
	Then we define
	\begin{align*}
		\bK^{(i)}_{\alpha} = \sbf^{(i)}_p\rbr{t;c_q;\boldsymbol{U}_{\alpha}^{(1)}(\bmu);\boldsymbol{U}_{\alpha}^{(2)}(\bnu)},\quad
		\bL^{(i)}_{\alpha} = \sbg^{(i)}_p\rbr{t;c_q;\boldsymbol{U}_{\alpha}^{(1)}(\bmu);\boldsymbol{U}_{\alpha}^{(2)}(\bnu)}.
	\end{align*}
	Here, we have defined the functions $ \sbf^{(1)}$ and $\sbf^{(2)}$ by
	\begin{align*}
		& \sbf^{(1)}_p(t;c_q;\sbv;\sbw)(k) 
		\\ &= -i\lam\sum_{k+k_1 = k_2 + k_3}\frac{1}{t^{p+1}}\int_0^{c_q t}\mcK_{2d}(s;k,k_1,k_2,k_3)s^p dsw_{k_1}v_{k_2}v_{k_3},\\ &
		\sbf^{(2)}_p(t;c_q;\sbv;\sbw)(k) 
		\\ &= i\lam\sum_{k+k_1 = k_2 + k_3}\frac{1}{t^{p+1}}\int_0^{c_q t}\mcK_{2d}(s;k_2,k_3,k,k_1)s^p dsv_{k_1}w_{k_2}w_{k_3},
	\end{align*}
	and likewise, we defined $ \sbg^{(1)}$ and $\sbg^{(2)}$ by
	\begin{align*}
		& \sbg^{(1)}_p\rbr{t;c_q;\sbv;\sbw}(k)
		\\ &= -i\kappa\sum_{k=k_1+k_2}\frac{1}{t^{p+\frac12}}\int_0^{c_q t}  P_d[e^{i \cdot (kk_2 + k_1k_2)}](s) s^p\circ dW_{k_2}(s)v_{k_1}\Phi_{k_2},\\
		& 	\sbg^{(2)}_p\rbr{t;c_q;\sbv;\sbw}(k) \\ &= i\kappa\sum_{k=k_1+k_2}\frac{1}{t^{p+\frac12}}\int_0^{c_q t}  P_d[e^{-i \cdot (kk_2 + k_1k_2)}](s) s^p\circ dW_{k_2}(s)w_{k_1}\Phi_{k_2}.
	\end{align*}
	With this in mind, for each Fourier component,
	\begin{align*}
		\bmu^{n+1}_k &= e^{- i t k^2}\bmu^n_k + \sum_{\alpha\in D} e^{-i t k^2}\rbr{t  b_{\alpha}^{(0)} \bK^{(1)}_{\alpha}+\sqrt{t} b_{\alpha}^{(1)}\bL^{(1)}_{\alpha}}, \\
		\bnu^{n+1}_k &= e^{i t k^2}\bnu^n_k + \sum_{\alpha\in D} e^{i t k^2}\rbr{t b_{\alpha}^{(0)}\bK^{(2)}_{\alpha}+\sqrt{t} b_{\alpha}^{(1)}\bL^{(2)}_{\alpha}}.
	\end{align*}
	Before proving Proposition \ref{prop:conservation-of-energy}, we introduce a technical lemma that will do much of the work for us. We recall a definition before this lemma.
	For $ s\in\mbbR_+ $, we define the following norm on sequences $ \boldsymbol{a} = \set{a_m}_{m\in \mathbb{Z}} $:
	\begin{equation*}
		\norm{\boldsymbol{a}}_{2,s}:=\rbr{\sum_{m\in\mbbZ}\langle m\rangle^{2s}\abs{a_m}^2}^{1/2},
	\end{equation*}
	where $ \inner{\cdot} = \sqrt{1+\abs{\cdot}^2} $. The space of $ l^2_s $-sequences is then defined by 
	\begin{equation*}
		l^2_s(\mathbb{Z}) : = \set{\boldsymbol{a}\in l^2\rbr{\mbbZ \backslash \set{0}}:\norm{\boldsymbol{a}}_{2,s}<\infty}.
	\end{equation*}
	\begin{lemma}\label{lemma:vwedgeGiszero}
		For any $ \sbv,\sbw\in l^2_s(\mbbZ),s>1 $, one has
		\begin{align*}
			&d\sbv \wedge d\rbr{\sbf^{(1)}_p\rbr{t;c_q;\sbv;\sbw}} + d\rbr{\sbf^{(2)}_p\rbr{t;c_q;\sbv;\sbw}} \wedge d\sbw \\
			&+ d\sbv \wedge d\rbr{\sbg^{(1)}_p\rbr{t;c_q;\sbv;\sbw}} + d\rbr{\sbg^{(2)}_p\rbr{t;c_q;\sbv;\sbw}} \wedge d\sbw = 0.
		\end{align*}
	\end{lemma}
	The proof of Lemma \ref{lemma:vwedgeGiszero} is an extension of the proof of \cite[Lemma 4.4]{M.S2023} to the newly introduced random variables. 
	\begin{proof}
		First, we consider the deterministic part of the scheme. We have,
		\begin{align*}
			& d\sbv \wedge d\rbr{\sbf^{(2)}_p(t;c_q;\sbv;\sbw)} \\
			&= i\lam\int_0^{c_qt}\sum_{k+k_1 = k_2 + k_3}\mcK_{2d}(s;k_2,k_3,k,k_1)(dv_k\wedge dv_{k_1})w_{k_2}w_{k_3}s^pds \\
			&+2i\lam\frac{1}{t^{p+1}}\int_0^{c_q t}\sum_{k+k_1 = k_2+k_3}\mcK_{2d}(s;k_2,k_3,k,k_1)(dv_k\wedge dw_{k_2})v_{k_1}w_{k_3}s^pds.
		\end{align*}
		We now observe that $ \mcK_{2d}(s;k_2,k_3,k,k_1) = \mcK_{2d}(s;k_2,k_3,k_1,k) $ while the wedge product is antisymmetric. Applying both to the first term shows that it cancels, and we obtain,
		\begin{align}
			& 	d\sbv \wedge d\rbr{\sbf^{(2)}_p(t;c_q;\sbv;\sbw)}\notag \\ 
			&= 2i\lam\frac{1}{t^{p+1}}\int_0^{c_q t}\sum_{k+k_1 = k_2+k_3}\mcK_{2d}(s;k_2,k_3,k,k_1)(dv_k\wedge dw_{k_2})v_{k_1}w_{k_3}s^pds\label{eq:sym_det_1}.
		\end{align}
		The same reasoning gives 
		\begin{align*}
			& d\rbr{\sbf^{(1)}_p(t;c_q;\sbv;\sbw)}\wedge d \sbw\notag \\ 
			&= -2i\lam\frac{1}{t^{p+1}}\int_0^{c_q t}\sum_{k+k_1 = k_2+k_3}\mcK_{2d}(s;k_2,k_3,k,k_1)\label{eq:sym_det_2}  (dv_k\wedge dw_{k_2})v_{k_1}w_{k_3}s^pd s.
		\end{align*}
		One obtains in the end
		\begin{equation*}
			d\sbv \wedge d\rbr{\sbf^{(2)}_p\rbr{t;c_q;\sbv;\sbw}} + d\rbr{\sbf^{(1)}_p\rbr{t;c_q;\sbv;\sbw}} \wedge d\sbw = 0.
		\end{equation*}
		Now we want to perform the same calculation on the stochastic parts. First we have 
		\begin{align*}
			&	d \sbv \wedge  d\rbr{\sbg^{(2)}_p(t;c_q;\sbv;\sbw)}  \\ & = i\kappa\sum_{k=k_1+k_2}\frac{1}{t^{p+\frac12}}\int_0^{c_q t}  P_d[e^{-i \cdot (kk_2 + k_1k_2)}](s) s^p\circ dW_{k_2}(s)(d v_k\wedge d w_{k_1})\Phi_{k_2}.
		\end{align*}
		To make this cancel, first note that the noise coefficients are real, $ W_{k_2}(s,\omega)\in\mbbR $, and perform the permutation $ k \leftrightarrow k_1 $. Next, note that this permutation enforces the map $ k_2 \mapsto -k_2 $, and recall that, by hypothesis, we have  $ \Phi^*_{k}=\Phi_{k} $, so we obtain
		\begin{align*}
			& d\sbv \wedge d\rbr{\sbg^{(2)}_p(t;c_q;\sbv;\sbw)} \\ & = i\kappa\sum_{k=k_1+k_2}\frac{1}{t^{p+\frac12}}\int_0^{c_q t}  P_d[e^{i \cdot (kk_2 + k_1k_2)}](s) s^p\circ dW_{k_2}(s)(d  v_{k_1}\wedge d w_{k})\Phi_{k_2},
		\end{align*}
		while on the other hand, we have
		\begin{align*}
			& d\rbr{\sbg^{(1)}_p(t;c_q;\sbv;\sbw)} \wedge d\sbw \\
			&= - i\kappa\sum_{k=k_1+k_2}\frac{1}{t^{p+\frac12}}\int_0^{c_q t}  P_d[e^{i \cdot (kk_2 + k_1k_2)}](s) s^p\circ dW_{k_2}(s)(d v_{k_1}\wedge dw_{k})\Phi_{k_2}.
		\end{align*}
		We obtain the following result:
		\begin{equation*}
			d\sbv \wedge d\rbr{\sbg^{(2)}_p(t;c_q;\sbv;\sbw)} + d\rbr{\sbg^{(1)}_p(t;c_q;\sbv;\sbw)} \wedge d\sbw = 0.
		\end{equation*}
	\end{proof}
	\begin{proof}[Proof of Proposition \ref{prop:conservation-of-energy}]
		We aim to show that the symplectic form \eqref{eq:symplectic-form} is preserved in each iteration of the Scheme \ref{scheme:stoch-runge-kutta}, i.e.
		\begin{equation*}
			d\bmu^{n+1}\wedge d\bnu^{n+1} =  d\bmu^{n} \wedge d\bnu^{n}.
		\end{equation*}
		We plug in our scheme, and we obtain
		\begin{align*}
			&	d\bmu^{n+1}\wedge d\bnu^{n+1}  = de^{-it(\cdot)^2}\bmu^n\wedge de^{i t(\cdot)^2}\bnu^n \\
			& + \sum_{\alpha\in D} d\rbr{e^{-it(\cdot)^2}\bmu^n}\wedge d\rbr{e^{i t(\cdot)^2}\rbr{t b_{\alpha}^{(0)}\bK^{(2)}_{\alpha}+\sqrt{t} b_{\alpha}^{(1)}\bL^{(2)}_{\alpha}}} \\
			& + \sum_{\alpha\in D} d\rbr{e^{-it (\cdot)^2}\rbr{t b_{\alpha}^{(0)}\bK^{(1)}_{\alpha}+\sqrt{t} b_{\alpha}^{(1)}\bL^{(1)}_{\alpha}}} \wedge d\rbr{e^{i t(\cdot)^2}\bnu^n} \\
			& + \sum_{\alpha, \tilde{\alpha}\in D} d\rbr{e^{-i t(\cdot)^2}\rbr{t b_{\alpha}^{(0)}\bK^{(1)}_{\alpha}+\sqrt{t} b_{\alpha}^{(1)}\bL^{(1)}_{\alpha}}}\wedge d\rbr{e^{i t (\cdot)^2}\rbr{t b_{\tilde{\alpha}}^{(0)}\bK^{(2)}_{\tilde{\alpha}}+\sqrt{t} b_{\tilde{\alpha}}^{(1)}\bL^{(2)}_{\tilde{\alpha}}}}.
		\end{align*}
		Now, by moving to Fourier space, we obtain
		\begin{align*}
			de^{-it(\cdot)^2}\bmu^n\wedge de^{i t(\cdot)^2}\bnu^n & = \sum_{k \in \mathbb{Z}} 	de^{-itk^2}\bmu^n_k\wedge de^{i tk^2}\bnu^n_k
			\\ & =  \sum_{k \in \mathbb{Z}} e^{-itk^2}e^{i tk^2}	d\bmu^n_k\wedge d\bnu^n_k = d\bmu^n\wedge d\bnu^n.
		\end{align*}
		Adding and subtracting some terms, one gets
		\begin{align*}
			& d\bmu^{n+1}\wedge d\bnu^{n+1} -  d\bmu^{n} \wedge d\bnu^{n} \\ &= \sum_{\alpha \in D} d\rbr{\bmu^n + t \sum_{\ti{\alpha} \in D}  a^{(0)}_{\alpha,\ti{\alpha}}\bK^{(1)}_{\ti{\alpha}} + \sqrt{t} \sum_{\ti{\alpha} \in D}  a^{(1)}_{\alpha,\ti{\alpha}}\bL^{(1)}_{\ti{\alpha}}} \wedge d\rbr{t b_{\alpha}^{(0)}\bK^{(2)}_{\alpha}+\sqrt{t} b_{\alpha}^{(1)}\bL^{(2)}_{\alpha}}\\
			& + \sum_{\alpha \in D}d\rbr{t b_{\alpha}^{(0)}\bK^{(1)}_{\alpha}+\sqrt{t} b_{\alpha}^{(1)}\bL^{(1)}_{\alpha}} \wedge d\rbr{\bnu^n + \sum_{\ti\alpha \in D} a^{(0)}_{\alpha,\ti{\alpha}}t\bK^{(2)}_{\ti{\alpha}} + a^{(1)}_{\alpha,\ti\alpha}\sqrt{t}\bL^{(2)}_{\ti{\alpha}}}\\
			&+\sum_{\alpha,\ti\alpha\in D} \rbr{b^{(0)}_{\alpha}b^{(0)}_{\ti\alpha} - b^{(0)}_{\alpha}a^{(0)}_{\alpha,\ti\alpha} - b^{(0)}_{\ti\alpha}a^{(0)}_{\ti\alpha,\alpha}}  t^2 d\bK^{(1)}_{\alpha} \wedge d\bK^{(2)}_{\alpha}\\
			& + \sum_{\alpha,\ti\alpha\in D} \rbr{b^{(0)}_{\alpha}b^{(1)}_{\ti\alpha} - b^{(0)}_{\alpha}a^{(1)}_{\alpha,\ti\alpha} - b^{(1)}_{\ti\alpha}a^{(0)}_{\ti\alpha,\alpha}}  t^{\frac{3}{2}} d\bK^{(1)}_{\alpha} \wedge d\bL^{(2)}_{\alpha}
			\\
			& + \sum_{\alpha,\ti\alpha\in D} \rbr{b^{(0)}_{\alpha}b^{(1)}_{\ti\alpha} - b^{(0)}_{\alpha}a^{(1)}_{\alpha,\ti\alpha} - b^{(1)}_{\ti\alpha}a^{(0)}_{\ti\alpha,\alpha}}  t^{\frac{3}{2}} d\bL^{(1)}_{\alpha} \wedge d\bK^{(2)}_{\alpha}
			\\ &  + \sum_{\alpha,\ti\alpha\in D} \rbr{b^{(1)}_{\alpha}b^{(1)}_{\ti\alpha} - b^{(1)}_{\alpha}a^{(1)}_{\alpha,\ti\alpha} - b^{(1)}_{\ti\alpha}a^{(1)}_{\ti\alpha,\alpha}}  t d\bL^{(1)}_{\alpha} \wedge d\bL^{(2)}_{\alpha}.
		\end{align*}
		We conclude by Hypothesis \ref{hyp:coefficient} and Lemma \ref{lemma:vwedgeGiszero} that  we have 
		\begin{equation*}
			d\bmu^{n+1}\wedge d\bnu^{n+1} -  d\bmu^{n} \wedge d\bnu^{n}  = 0.
		\end{equation*}
	\end{proof}
	\section{Convergence Analysis of the Midpoint Rule}\label{sec::6}
	In the previous sections, we have developed the theory for a class of $ \mcO\rbr{\sqrt t} $ resonance-based Runge-Kutta schemes for the non-linear Schr\"odinger equation. We will now look closely at the convergence of one particular scheme, specifically the one introduced in Example \ref{ex:stoch-res-midpoint}. Recall that in this example, we derived the stochastic resonance-based midpoint rule, which is written for initial conditions $ u_0 $ as follows:
	\begin{equation}\label{scheme:midpoint-rule}
		S^{n+1}_t \rbr{u_0} = u_{n+1} = e^{it\pa_x^2}u_n + e^{it\pa_x^2}\rbr{t\mathfrak{F}_{0}\rbr{t;1;u_{n+\frac12}} + \sqrt{t}\mathfrak{P}_{0}\rbr{t;1;u_{n+\frac12}}}
	\end{equation}
	where
	\begin{equation*}
		u_{n+\frac12} = \frac{u_n + e^{-it\pa_x^2}u_{n+1}}{2}.
	\end{equation*}
	and the random variables in the scheme are 
	\begin{equation*}
		\sqrt{t}W_n(x) = \mcF^{-1}\rbr{\sum_{k\in\mbbZ}e^{ikx}\int_{t_n}^{t_{n+1}}\circ dW(s)}.
	\end{equation*}
		
	For conservation of mass and energy we require real-valued noise $ W_n(x) $. This is equivalent to symmetry of frequencies in Fourier space, i.e, the condition that the Fourier coefficients satisfy 
	\begin{equation}
		W_n(-k) = \overline{W_n(k)} = W_n(k).
	\end{equation}
	Thus, when the noise is real, the modulus $ |W_n(x)| $ is simply the absolute value of the real-valued process $ W_n(x) $.
	If $ W_n(x) $ is complex, $ W_n(x)\overline{W_n(x)} = |W_n(x)|^2 $, hence, essentially the same notation extends to a complex valued process for which our scheme still converges, albeit without the conservation properties.
	
	To perform fixed-point iterations, we will introduce an auxiliary function with fixed randomness. The stochastic discretisation map will be altered to
	\begin{equation}\label{eq:aux-func-stoch}
		\begin{aligned}
			\mcP^*_0(t;1;v,X) & = -i  \kappa\sum_{k\in\mbbZ}e^{ixk}\frac{1}{\sqrt{t}}\sum_{k = k_1 + k_2}v_{k_1}\Phi_{k_2}X_{k_2},
		\end{aligned}
	\end{equation}
	where $ X $ is a fixed random variable. Now we define the auxiliary process
	\begin{align*}
		\mathfrak{S}_t(w)&:=e^{it\pa_x^2}u_n + te^{it\pa_x^2}\mfF_0\rbr{t;1;u_{n+\frac12}^w} + \sqrt{t}e^{it\pa_x^2}\mfP^*_0\rbr{t;1;u_{n+\frac12}^w,W_n},
	\end{align*} 
	where $ W_n(t) $ is a centred Gaussian vector with variance $ t $, inherited from the $ n^{\text{th}} $ iteration of our scheme, $ S^n_t(u_0) $. To simplify the notation, we have defined: 
	\begin{equation*}
		u_{n+\frac12}^w =  \frac{u_n + e^{-it\pa_x^2}w}{2}.
	\end{equation*}
	Note that iterations of $ \mfS_t(\cdot) $ are not random, because the random variable is fixed from the beginning. In the next subsection, we will demonstrate that each iteration of Scheme \eqref{scheme:midpoint-rule} has an almost sure fixed point, provided the timestep is sufficiently small. 
	\subsection{Existence and Stability}
	Because the methods discussed in this paper are implicit, we must employ a fixed-point iteration to solve the method at each time step; for this reason, we aim to demonstrate that such a process indeed converges to the correct solution. This can be captured in the following statement:
	\begin{theorem}[Convergence of Fixed Point Iterations]\label{thm:fixedpoint-convergence}
		Let $ R_n>0 $ and $ \alpha >1 $. Then, for any $ u_n,v,w\in B_{R_n}(H^\alpha) = \set{u\in H^\alpha : \norm{u}_{H^\alpha}< R_n} $ we have that $ \mfS_t^j(u_n) $ converges to $ u_{n+1} $ almost surely, i.e.
		\begin{equation*}
			\mbbP\rbr{\lim_{j\rightarrow\infty}\mfS_t^j(u_n) = u_{n+1}}=1.
		\end{equation*}
	\end{theorem}
	\begin{proof}
		We want to show that for sufficiently small, but positive $ t $, that the iterates of $ \mfS_t(u_n) $ form a convergent sequence. We have  
		\begin{align*}
			&\norm{\mfF_0(t;1;u_{n+\frac12}^w) - \mfF_0(t;1;u_{n+\frac12}^v)}^2_{H^\alpha} \\
			&=t^2\frac{\abs\mu^2}{2}\sum_{k\in\mbbZ} \langle k \rangle^{2\alpha}\bigg|\sum_{k+k_1 = k_2 + k_3}\sqbr{\frac{e^{-2itkk_1} - 1}{-2itkk_1} + \frac{e^{2itk_2k_3} - 1}{-2itk_2k_3} - 1} \\
			&\bigg(\rbr{\bar{u}^n_{k_1} + \bar{w}_{k_1}}\rbr{{u}^n_{k_2} + {w}_{k_2}}\rbr{{u}^n_{k_3} + {w}_{k_3}} - \rbr{{\bar{u}}^n_{k_1} + {\bar{v}}_{k_1}}\rbr{{u}^n_{k_2} + {v}_{k_2}}\rbr{{u}^n_{k_3} + {v}_{k_3}}\bigg)\bigg|^2.
		\end{align*}
		By applying the following factorisation to the non-linearity,
		\begin{equation} \label{factorisation_abc}
			\begin{aligned}
				abc - def &= abc - bcd + bcd -edc + cde -efd \\
				& = (a-d)bc + (b-e)dc + (c-f)de, 
			\end{aligned}
		\end{equation}
		and observing that the exponential functions in $ k,k_1,k_2 $ and $ k_3 $ are uniformly bounded, we have,  by the triangle inequality, that 
		\begin{equation*}
			t\norm{\mfF_0(t;1;u^v_{n+\frac12}) - \mfF_0(t;1;u^w_{n+\frac12})}_{H^\alpha}\leq tC_{R_n}\norm{v-w}_{H^\alpha}.
		\end{equation*}
		Then, for the stochastic term, we have 
		\begin{equation*}
			\sqrt{t}\norm{\mfP^*_0(t;1;u^v_{n+\frac12},X) - \mfP^*_0(t;1;u^w_{n+\frac12},X)}_{H^\alpha} \leq \sqrt{t}\Phi \abs{X}\norm{v-w}_{H^\alpha}.
		\end{equation*} 
		This implies that for $ v,w\in B_{R_n}(H^\alpha) $, we have 
		\begin{equation*}
			\norm{\mfS_t(w) - \mfS_t(v)}_{H^\alpha} \leq \rbr{C_{R_n}t + C_{\Phi,X}\sqrt{t}}\norm{v-w}_{H^\alpha}.
		\end{equation*}
		It follows that
		\begin{equation*}
			\norm{\mfS_t(e^{it\pa_x^2}w) - \mfS_t(e^{it\pa_x^2}v)}_{H^\alpha}\leq\rbr{C_{R_n}t + \sqrt{t}C_{\Phi,W_n}}\norm{v-w}_{H^\alpha},
		\end{equation*}
		where $ R_n $ is the radius of the ball containing $ u_n $, and $ W_n $ is the random variable introduced by $ S^n_t(u_{0}) $. Now, selecting $ t_n $ such that $ t_nC_{R_n} + \sqrt{t_n}C_{\Phi,W_n}<1 $, i.e.
		\begin{equation}\label{eq::bound-on-tn}
			0<t_n<\frac{C_{\Phi ,W_n}^2+2 C_{R_n}}{2 C_{R_n}^2}-\frac{1}{2} \sqrt{\frac{4 C_{R_n} C_{\Phi ,W_n}^2+C_{\Phi ,W_n}^4}{C_{R_n}^4}},
		\end{equation}
		we apply the Banach fixed point theorem to prove the result. Note that this holds almost surely, as the random variable is infinite with negligible probability, but it is nonetheless possible.
		\begin{remark}
			As a result of the condition \eqref{eq::bound-on-tn}, the scheme can not be adapted and well-posed with probability $ 1 $. To be well-posed, we require the fixed point iterations to converge and to make this probability $ 1 $ we need to select $ t_n $ based on the random variables in the $ n^{\text{th}} $ iteration, which is only $ \mcF_{t_{n+1}}- $adapted.
		\end{remark}
	\end{proof}
	\begin{remark}[Problem concerning fixed point convergence] \, \newline
		We note that the bound \eqref{eq::bound-on-tn} is very restrictive. In fact, it is clear that due to the randomness in equation \eqref{eq::snlse}, for any fixed time step $ t $, the scheme will be ill-posed with non-negligible probability. This problem is observed in \cite{B.D2004}, where convergence of a discretisation is proved as $ t\rightarrow 0$. The scheme in \cite{B.D2004} has been successfully applied in simulations of the SNLSE in \cite{D.D2002}. To prove weak convergence of similar schemes, Milstein, Repin, and Tretyakov \cite{MRT2002} use a truncated noise; this approach does not allow strong convergence, however, as the error between the truncated variable and the randomness driving \eqref{eq::snlse} can be arbitrarily large. 	
	\end{remark}
	In the next theorem, we will show that the stability of the scheme requires the same condition as the fixed-point convergence. 
	\begin{theorem}[Stability]\label{thm:stability}
		Fix $ \alpha > 1 $ and let $ S_{t}^n(u) = u_n $. For $ 0\leq n \leq N $, define $ R_n = \min\set{R\in\mbbN : v_n,w_n\in B_R(H^\alpha)}$ and $ R^* = \max\set{R_n<\infty: n\leq N}$. Define the stopping time $ \tau = tn^* $, which occurs if $ R_{n^*}=\infty $. Let $ T=Nt $ be the length of the simulation, and finally, let $ m = \min\set{N,n^*} $. Then we have 
		\begin{equation*}
			\lim_{\norm{v - w}_{H^\alpha}\rightarrow 0}\sup_{0\leq n\leq m}\mbbP \rbr{\norm{S_t^n(v) - S^n_t(w)}_{H^\alpha}\geq K} =0.
		\end{equation*}
	\end{theorem}
	\begin{proof}[Proof of Theorem \ref{thm:stability}]
		By the hypothesis of the theorem, we have $ v,w\in B_{R_0}(H^\alpha) $. Setting $ \ti{v} = \frac{v+e^{-it\pa_x^2}S_t(v)}{2}$ and $ \ti{w} = \frac{w+e^{-it\pa_x^2}S_t(w)}{2}$, we find that
		\begin{align*}
			\norm{S_t(v) - S_t(w)}_{H^\alpha} &\leq \norm{e^{it\pa_x^2}(v-w)}_{H^\alpha} \\
			&+t\norm{e^{it\pa_x^2}\mfF_0\rbr{c_q;1;\ti v} - e^{it\pa_x^2}\mfF_{0}\rbr{c_q;1;\ti w}}_{H^\alpha} \\
			&+ \sqrt{t}\norm{e^{it\pa_x^2}\mfP_0\rbr{c_q;1;\ti v} - e^{it\pa_x^2}\mfP_{0}\rbr{c_q;1;\ti w}}_{H^\alpha}.
		\end{align*}
		Clearly the first term is bounded by a constant $ C_{R_0} $ since $ e^{it\pa_x^2} $ is unitary on $ H^\alpha $. For the second term, we follow the argument laid out in \cite{M.S2023}, which we reproduce here for convenience. First, we note that 
		\begin{align*}
			\norm{\tilde{v}-\tilde{w}}_{H^\alpha} &  = 
			\norm{  \frac{v+e^{-it\pa_x^2}S_t(v)}{2} - \frac{w+e^{-it\pa_x^2}S_t(w)}{2}}_{H^\alpha}
			\\ & \leq \frac{1}{2} \left( \norm{ v-w }_{H^\alpha} + \norm{ S_t(v)-S_t(w) }_{H^\alpha}  \right).
		\end{align*}
		Then since $ w_n,v_n,S_t(w_n), S_t(v_n)\in B_{R^*}(H^\alpha) $, we have 
		\begin{align*}
			t\norm{e^{it\pa_x^2}\mfF\rbr{t;1; v} - e^{it\pa_x^2}\mfF\rbr{t;1;w}}_{H^\alpha} \leq tC_{R^*}\norm{\rbr{v-w} + S_t(v) - S_t(w)}_{H^\alpha}.
		\end{align*}
		For the stochastic term, we use the estimate 
		\begin{align*}	
			\,	&	\sqrt{t}\norm{\mfP_0(t;1;v_{n+\frac12}) - \mfP_0(t;1;w_{n+\frac12})}_{H^\alpha}\leq
			\\ & \sqrt{t} C_{\Phi,W_n}\rbr{\norm{v_n - w_n}_{H^\alpha} + \norm{S_t(v_n) - S_t(w_n)}_{H^\alpha}}.
		\end{align*}
		Putting these together, we obtain
		\begin{align*}
			\norm{S_t(v) - S_t(w)}_{H^\alpha} &\leq \norm{v-w}_{H^\alpha} + t\frac{C_{R^*}}{2} \left( \norm{ v-w }_{H^\alpha} + \norm{ S_t(v)-S_t(w) }_{H^\alpha}  \right) \\
			&\quad + \sqrt{t}\frac{|W|}{2} \left( \norm{ v-w }_{H^\alpha} + \norm{ S_t(v)-S_t(w) }_{H^\alpha}  \right),
		\end{align*}
		where $ W $ is a standard Gaussian random variable. This implies the inequality
		\begin{equation}\label{eq:stability_ineq}
			\begin{aligned}
				\norm{S_t(v) - S_t(w)}_{H^\alpha} &\leq\frac{1 + \frac{tC_{R} + \sqrt{t}|W|}{2}}{1-\frac{tC_{R} + \sqrt{t}|W|}{2}}\norm{v-w}_{H^\alpha}.
			\end{aligned}
		\end{equation}
		\begin{remark}
			Note that the function $ \frac{1+x}{1-x} $ blows up at $ x=1 $ and is negative for values outside of $ [0,1) $. Since the right side is positive, the inequality in \eqref{eq:stability_ineq} makes sense only if \[tC_{R_n} + \sqrt{t}\abs{W_n} < 1.\]
			This is the precise condition required on $ t $ for the fixed point convergence.
		\end{remark}
		Next, we use the fact that 
		\begin{equation*}
			\frac{1+x/2}{1-x/2}\leq e^{3x/2},\quad \text{for }x<1,
		\end{equation*}
		and note that $ e^{\frac32 \rbr{tC_{R_n} + \sqrt{t}\abs{W_n}}} > 1 $. Iterating the scheme and applying inequality \eqref{eq:stability_ineq} at each step leads to the following inequality
		\begin{equation*}
			\norm{S_t^n(v) - S_t^n(w)}_{H^\alpha} \leq \exp\rbr{X_n}\norm{v-w}_{H^\alpha},
		\end{equation*}
		where we define $ X_n = \sum_{m=1}^n \rbr{tC_{R_m} + \sqrt{t}\abs{W_m}} $. It follows that $ e^{X_n} $ has increasing expectation and is hence a sub-martingale in terms of $ n $, and thus we may apply Markov's inequality and Doob's maximal inequality to obtain
		\begin{equation*}
			\mbbP\rbr{\sup_{n\leq m}\norm{S_t^n(v) - S_t^n(w)}_{H^\alpha}\geq K} \leq \frac{2\mbbE\sqbr{e^{X_m}}\norm{v-w}_{H^\alpha}}{K} \quad a.s.
		\end{equation*}
	\end{proof}
	\subsection{Local Error}
	We have seen that the scheme is unstable with a non-negligible probability; however, we can prove the following result about local error. 
	\begin{theorem}[Local Error]\label{prop:local-error}
		Denote by $ E_t(u(t_n)) $ the exact solution at time $ t_n + t $ from the initial datum $ u(t_n) $. Let $ \alpha>1 $ and choose $ R>0 $ such that \[ \sup_{s\in[0,t]}\norm{E_s(u(t_n))}_{H^{\alpha}} \leq R, \] and that $ u(t_n)\in B_{R}(H^\alpha) $, then
		\begin{equation*}
			\norm{S_t(u(t_n)) - E_t(u(t_n))}_{H^\alpha}\leq C_1t^{3/2},
		\end{equation*}
		where $ C_1 > 0$ is a random constant depending on $ \Phi $, $ \alpha $, $ W_1 $ and $ W_2^2 $ for 
		\begin{equation*}
			W_1 = \frac{1}{\sqrt{t}}\int_{0}^{t}\circ dW(s),\quad W_2 = \frac{1}{t}\int_{0}^{t}\int_{0}^{r}\circ dW(s)\circ dW(r).
		\end{equation*}
	\end{theorem}
	\begin{proof}

	Note that by the Burkholder-Davis-Gundy inequality and the assumption that $\Phi$ is Hilbert-Schmidt, one can show that $\Phi W(t,x)$ is almost surely bounded in $H^s$ for each $t$. Specifically, for any $T > 0$:
	$$
	\mathbb{E}\left[\sup_{0 \leq t \leq T} \|\Phi W(t, \cdot)\|_{H^s}^2\right] \lesssim T \|\Phi\|_{\text{HS}(L^2, H^s)}^2.
	$$
	This implies that $\Phi W(t,x)$ is almost surely in $H^s$ for each $t$, which will hold for the remainder terms in the series expansions of our operators. Next, we introduce the shorthand notations 
		\begin{equation*}
			w_{t_n+\frac12}(s) =e^{is\pa_x^2}\frac{u(t_n) + e^{-it\pa_x^2}S_t(u(t_n))}{2}, \quad  \tilde{w}_{t_n+\frac12} =\frac{u(t_n) + e^{-it\pa_x^2}S_t(u(t_n))}{2}.
		\end{equation*} 
		Initially, we will look at the difference between the scheme in \ref{scheme:stoch-runge-kutta} and the Duhamel iteration of $ w_{t_n+\frac12}(s) $, meaning we want to bound the $ H^\alpha $  norm of the following terms.
		\begin{align*}
			\mcE_1(u(t_n)) & = \int_0^t e^{i(t-s)\pa_x^2} |u(t_n+s)|^2 \bar{u}(t_n+s) ds  - t\mfF_{0}\rbr{t;1;\ti{w}_{n+\frac12}},
			\\
			\mcE_2(u(t_n)) & = \int_0^t e^{i(t-s)\pa_x^2} u(t_n+s)\circ \Phi dW(s) - \sqrt{t}\mfP_{0}\rbr{t;1;\ti{w}_{n+\frac12}}.
		\end{align*}
		The primary challenge in the analysis arises from the presence of midpoint approximations in the maps $\mfF$ and $\mfP$, while the corresponding integrals are formulated using the exact solution. To overcome this, we employ the following representation:
		\begin{align*}
			& 	\norm{S_t(u(t_n)) - E_t(u(t_n))}_{H_\alpha}  \\ &= \norm{S_t(u(t_n)) -  E_t(w_{t_n+\frac12}) + E_t(w_{t_n+\frac12})- E_t(u(t_n))}_{H^\alpha} \\
			&\leq \norm{S_t(u(t_n)) -  E_t(w_{t_n+\frac12})}_{H^\alpha} + \norm{E_t(w_{t_n+\frac12}) - E_t(u(t_n))}_{H^\alpha},
		\end{align*}
		where we define
		\begin{align*}
			E_t(w_{t_n+\frac12}(s)) &= e^{i(t_n + t)\pa_x^2}u(t_n) - i\kappa e^{it\pa_x^2}\int_0^t e^{-ts\pa_x^2}\abs{w_{n+\frac12}(s)}^2w_{n+\frac12}(s)ds \\
			&\quad -i \lam e^{it\pa_x^2}\int_0^t e^{-is\pa_x^2} w_{n+\frac12}(s)\Phi\circ dW(s).
		\end{align*}
		We can now easily bound the first term in the sum 
		\begin{align*}
			& E_t(w_{t_n+\frac12}(s)) - S_t(u(t_n)) \\ &=\int_0^t e^{-i(t_n+s)\pa_x^2} |w_{t_n+\frac12}(s)|^2 w_{t_n+\frac12}(s) ds - t\mfF_{0}\rbr{t;1;\tilde{w}_{t_n+\frac12}} \\
			&+ \int_0^t e^{-i(t_n+s)\pa_x^2} w_{t_n+\frac12}(s)\circ \Phi dW(s) - \sqrt{t}\mfP_{0}\rbr{t;1;\tilde{w}_{t_n+\frac12}} \\
			&= \mcE_{1,k}(w_{t_n+\frac12}) + \mcE_{2,k}(w_{t_n+\frac12}).
		\end{align*}
		Applying the Fourier transform to the first term gives rise to the following integral 
		\begin{align*}
			\mcE_{1,k}(w_{t_n+\frac12})&= \sum_{k=-k_1+k_2+k_3} \overline{\tilde{w}^n_{k_1}}\tilde{w}^n_{k_2}\tilde{w}^n_{k_3}\\ & \int_0^t e^{-2iskk_1+2isk_2k_3} - \sqbr{e^{-2iskk_1} + e^{2isk_1k_3}-1}ds,
		\end{align*}
		where $ \tilde{w}^n_k $ (resp. $\mcE_{1,k}(w_{t_n+\frac12})$) is the $k^{\text{th}}$ Fourier coefficient of $ \tilde{w}_{t_n+\frac12}  $ (resp. $\mcE_{1}(w_{t_n+\frac12})$). Then
		\begin{equation*}
			\mcE_{1,k}(w_{t_n+\frac12}) =  \mathcal{O}\rbr{t^{2} \sum_{k=-k_1+k_2+k_3} (|k k_1 k_2 k_3|)^{1/2}\overline{\tilde{w}_{k_1}}\tilde{w}_{k_2}\tilde{w}_{k_3}}, \quad a.s.
		\end{equation*}
		Here, we have used the error estimation on the kernel approximation given in \eqref{bound_kernel}.
		For the stochastic integral, we obtain the estimate
		\begin{align*}
			\mcE_{2,k}(w_{t_n+\frac12}) & = \sum_{k=k_1 + k_2} \tilde{w}_{k_1}\Phi_{k_2} \int_0^t \rbr{e^{is\rbr{k_2^2 + k_1k_2}} - 1 }\circ dW_{k_2}(s).
			\\ & = \mathcal{O}\rbr{t^{3/2} \sum_{k=k_1+k_2}(k_2^2 + k_1k_2)w_{k_1}\Phi_{k_2}W_{k_2}},\quad a.s,
		\end{align*}
		where the error is a consequence of the Taylor approximation.
		Combining these estimates, we almost-surely get,
		\begin{equation*}
			\norm{E_t(w_{t_n+\frac12}(t_n)) - S_t\rbr{u(t_n)}}_{H^\alpha}  = \mathcal{O}\rbr{t^{\frac{3}{2}} \text{Tr}(\lap\Phi W_1) \norm{u}_{H^\alpha}}.
		\end{equation*}
		We now need to bound the term 
		\begin{equation*}
			\norm{E_t(w_{t_n+\frac12}) - E_t\rbr{u(t_n)}}_{H^\alpha}.
		\end{equation*}
		The first step is to rewrite $ E_t\rbr{u(s)} $ in terms of a midpoint approximation of $ u(s) $, allowing us to write 
		\begin{align*}
			&\norm{E_t\rbr{w_{t_n+\frac12}} - E_t\rbr{u(t_n)}}_{H^\alpha}\\
			&\leq \bigg|\bigg|\mu \int_{0}^{t} e^{i(t-s)\pa_x^2} \bigg(|u(t_n+s)|^2 u(t_n+s) -|w_{t_n+\frac12}(s)|^2 w_{t_n+\frac12}(s)\bigg)ds\bigg|\bigg|_{H^\alpha} \\
			&\quad+\norm{\kappa \int_{0}^{t} e^{i(t-s)\pa_x^2}\rbr{u(t_n+s) - w_{t_n+\frac12}(s)} \circ \Phi d W(s)}_{H^\alpha}.
		\end{align*}
		For the first term we can use the factorisation \eqref{factorisation_abc} presented in Lemma \ref{thm:fixedpoint-convergence} for the non-linearity, almost-surely giving
		\begin{equation*}
			\bigg|\bigg|\mu \int_{0}^{t} e^{i(t-s)\pa_x^2} \bigg(|u(t_n+s)|^2 u(t_n+s) -|w_{t_n+\frac12}(s)|^2 w_{t_n+\frac12}(s)\bigg)ds\bigg|\bigg|_{H^\alpha} = \mcO(t^{3/2}).
		\end{equation*} 
		For the stochastic term, we first rewrite $ u(t_n+s) $ using the following identities:
		\begin{align*}
			u(t_n+s) &= e^{i s\pa_x^2} u(t_n)   - i \kappa e^{i s\pa_x^2} \int_0^{s} e^{-i r\pa_x^2}u(t_n+r) \circ\Phi d W(r)  + \mathcal{O}(s)
			\\
			u(t_n+s) &=  e^{i (s-t)\pa_x^2} u(t_n+t)  - i \kappa e^{i s\pa_x^2} \int_t^s e^{-i r\pa_x^2}u(t_n+r) \circ \Phi d W(r) + \mathcal{O}(s).
		\end{align*}
		Combining gives
		\begin{align*}
			u(t_n+s) &=  w_{t_n+\frac12}(s)  - i \kappa  \bigg(\frac12\int_t^s e^{i(s-r)\pa_x^2}u(t_n+r) \circ \Phi d W(r) \\
			&\quad + \frac12\int_0^se^{i(s-r)\pa_x^2}u(t_n+r) \circ \Phi d W(r)\bigg) + \mathcal{O}(s) \quad a.s.
		\end{align*}
		Now, one almost-surely has
		\begin{align*}
			& i \kappa \int_{0}^{t}  e^{i(t-s)\pa_x^2}\rbr{u(t_n+s) - w_{t_n+\frac12}(s)} \circ \Phi d W(s) \\
			&= - \kappa^2\int_{0}^{t} e^{i(t-s)\pa_x^2}\bigg[  \bigg(\frac12\int_t^s e^{i(s-r)\pa_x^2}u(t_n+r) \circ \Phi d W(r) \\
			&\quad + \frac12\int_0^se^{i(s-r)\pa_x^2}u(t_n+r) \circ \Phi d W(r)\bigg) + \mathcal{O}(s)\bigg]\circ \Phi d W(s) \\
			& = \mcO\rbr{t^{2} \text{Tr}(\lap\Phi W_2)^2 \norm{u_n}_{H^\alpha}} +\mcO(t^{\frac{3}{2}})
		\end{align*}
		where the order $ t^{\frac{3}{2}} $ term is the truncation error of the Duhamel formula and the order $ t^2 $ term is derived precisely as in the calculation of \eqref{eq:order-of-iterate}. In the end, we obtained
		\begin{equation*}
			\norm{E_t(w_{t_n+\frac12}) - E_t\rbr{u(t_n)}}_{H^\alpha} =  \mcO(t^{\frac{3}{2}} \text{Tr}(\lap\Phi W_1) \text{Tr}(\lap\Phi W_2)^2 \norm{u_n}_{H^\alpha}),
		\end{equation*}
		or equivalently
		\begin{equation*}
			\norm{E_t(w_{t_n+\frac12}) - E_t\rbr{u(t_n)}}_{H^\alpha} \leq t^{\frac{3}{2}} C_{R, \Phi, W_1, W_1^2 }, \quad a.s,
		\end{equation*}
		which allows us to conclude.
	\end{proof}
	\section{Numerical Experiments}
	We have implemented the numerical scheme \ref{scheme:midpoint-rule} and demonstrate strong and pathwise convergence. For comparison we also implemented the stochastic Crank-Nicolson method studied in \cite{D.D2002}. 
	
	We find that the Crank-Nicolson Scheme has generally better mass preservation in the stochastic case but only works for modestly discretised spatial mesh and small timesteps. On the other hand, the Resonance based method has improved conservation properties for a more refined spatial mesh, and good stability for larger timesteps. For example, our scheme converges for spatial mesh of size $ 2^{-13} $ (8192 Fourier modes) with a timestep as large as $ 10^{-2} $. With the same timestep, however, the Crank-Nicolson scheme can converge with at maximum a mesh of the size $ 2^{-3} $ (8 Fourier modes).  This makes our scheme much better suited to the study of long time scales on refined spatial meshes. This scenario is applicable to many modern topics in dispersive (S)PDE such as energy cascading and ergodicity of solutions. 
	
	We note that in practice both schemes are stable and convergent, which in theory is not true as the schemes are almost surely unstable after a large enough period of time. On the other hand neither scheme are symplectic to machine precision, which \emph{is} true in theory. In other words, we find substantial mismatch between theory and practice which are likely the result of finite computational resources.  
	\subsection{Pathwise Convergence}
	The plots in this section show a single random path averaged over random initial data in $ H^s(\mbbT) $ with $ 64 $ frequencies. In figure \ref{fig:convergence_pathwise} we see that at low regularity, the order of convergence of the Crank-Nicolson is lowered significantly, while the resonance-based scheme maintains convergence of order $t$. We demonstrate that the rate of convergence is comparable between the two schemes when the regularity of the initial data is sufficient.

\begin{figure}
	\centering
	\begin{tikzpicture}
		\begin{groupplot}[
			group style = {
				group size=2 by 1,
				horizontal sep=2cm, 
			},
			width=0.5\textwidth,
			height=0.5\textwidth,
			grid=both,
			minor tick num=1,
			major grid style={lightgray},
			minor grid style={lightgray!25},
			xlabel={Timestep (t)},
			ylabel={$L^2$ Error},
			]
			\nextgroupplot[
			title={ Pathwise Convergence in $H^1$},
			title style={yshift=+.2cm}, 
			xmin=5e-6, xmax=1e-3,
			ymin=5e-6, ymax=1,
			xmode = log,
			ymode = log,
			]   
			\addplot coordinates {
				(1.000000e-3, 4.309876e-03)
				(5.000000e-4, 1.362383e-03)
				(1.000000e-4, 1.667492e-04)
				(5.000000e-5, 8.045209e-05)
				(1.000000e-5, 1.511307e-05)
				(5.000000e-6, 6.850583e-06)
			};
			
			\addplot coordinates {
				(1.000000e-3, 1.809680e-01)
				(5.000000e-4, 1.201644e-01)
				(1.000000e-4, 6.435442e-03)
				(5.000000e-5, 1.869023e-03)
				(1.000000e-5, 1.551735e-04)
				(5.000000e-6, 5.922330e-05)
			};
			
			\addplot[orange, dotted, thick, domain=1e-6:1e-3, samples=100] {10000000*x*x};

			\addplot[gray, dashed, thick, domain=1e-6:1e-3, samples=100] {7*x};
			
			\nextgroupplot[
			title={ Pathwise Convergence in $ H^{10}$},
			title style={yshift=+0.2cm}, 
			ylabel={$L^2$ Error},
			xlabel={Timestep (t)},
			xmin=5e-6, xmax=1e-3,
			ymin=5e-7, ymax=1e-3,
			xmode = log,
			ymode = log,
			]   
			\addplot coordinates {
				(1.000000e-3, 7.915679e-04)
				(5.000000e-4, 3.538177e-04)
				(1.000000e-4, 7.468616e-05)
				(5.000000e-5, 3.680867e-05)
				(1.000000e-5, 6.654243e-06)
				(5.000000e-6, 3.093619e-06)
			};
			
			\addplot coordinates {
				(1.000000e-3, 6.736100e-04)
				(5.000000e-4, 2.945556e-04)
				(1.000000e-4, 6.249015e-05)
				(5.000000e-5, 3.096632e-05)
				(1.000000e-5, 5.581342e-06)
				(5.000000e-6, 2.601243e-06)
			};
			
			\addplot[gray, dashed, thick, domain=1e-6:1e-3, samples=100] {x};
		\end{groupplot}
	\end{tikzpicture}
	\caption{With a spatial mesh of 64 frequencies we average over $ 10 $ instances of random initial data in $H^1(\mathbb{T})$ and $ H^{10}(\mathbb{T})$. We use a single path of multiplicative real noise.}
	\label{fig:convergence_pathwise}
\end{figure}
	\begin{figure}
	\centering
	\begin{tikzpicture}
		\begin{groupplot}[
			group style = {
				group size=2 by 1,
				horizontal sep=2cm, 
			},
			width=0.5\textwidth,
			height=0.5\textwidth,
			grid=both,
			minor tick num=1,
			major grid style={lightgray},
			minor grid style={lightgray!25},
			xlabel={Timestep (t)},
			ylabel={$L^2$ Error},
			legend style={
				at={(1.1,-0.3)},  
				anchor=north,
				legend columns=1,
				draw=none,
				font=\small,
			},
			]
			\nextgroupplot[
			title={Strong Convergence in $H^1$},
			title style={yshift=+0.2cm}, 
			xmin=5e-6, xmax=1e-3,
			ymin=1e-7, ymax=5e-1,
			xmode = log,
			ymode = log,
			]   
			\addplot coordinates {
				(1.000000e-3, 9.468936e-04)
				(5.000000e-4, 3.334848e-04)
				(1.000000e-4, 5.013668e-05)
				(5.000000e-5, 2.440155e-05)
				(1.000000e-5, 4.664105e-06)
				(5.000000e-6, 2.153750e-06)
			};
			\addlegendentry{Resonance Midpoint}
			
			\addplot coordinates {
				(1.000000e-3, 1.126137e-01)
				(5.000000e-4, 3.379715e-02)
				(1.000000e-4, 1.572000e-03)
				(5.000000e-5, 4.557991e-04)
				(1.000000e-5, 3.777176e-05)
				(5.000000e-6, 1.441330e-05)
			};
			\addlegendentry{Crank-Nicolson}
			
			\addplot[orange, dotted, thick, domain=1e-6:1e-3, samples=100] {1000000*x*x};
			\addlegendentry{Square: $y = x^2$}
			
			\addplot[gray, dashed, thick, domain=1e-6:1e-3, samples=100] {x/5};
			\addlegendentry{Linear: $y =  x$}
			
			\nextgroupplot[
			title={Strong Convergence in $H^{10}$},
			title style={yshift=+0.2cm}, 
			ylabel={$L^2$ Error},
			xlabel={Timestep (t)},
			xmin=5e-6, xmax=1e-3,
			ymin=5e-7, ymax=1e-3,
			xmode = log,
			ymode = log,
			]   
			\addplot[
			color=blue,
			mark=*,
			]
			coordinates {
				(1e-3, 2.087086e-04)
				(5e-4, 1.098482e-04)
				(1e-4, 2.010688e-05)
				(5e-5, 9.092293e-06)
				(1e-5, 1.838745e-06)
				(5e-6, 8.811336e-07)
			};
			
			\addplot[
			color=red,
			mark=square*,
			]
			coordinates {
				(1e-3, 1.972572e-04)
				(5e-4, 1.063453e-04)
				(1e-4, 1.961631e-05)
				(5e-5, 8.860516e-06)
				(1e-5,1.760371e-06)
				(5e-6, 8.395767e-07)
			};
			
			\addplot[
			color=gray,
			dashed,
			thick,
			domain=5e-6:1e-3,
			samples=100,
			]
			{x};
		\end{groupplot}
	\end{tikzpicture}
	\caption{With a spatial mesh of 64 frequencies we average over $ 33 $ instances of random initial data in $H^1(\mathbb{T})$ and $H^{10}(\mathbb{T})$. We use $ 33 $ paths of multiplicative real noise.}
	\label{fig:strong_convergence}
\end{figure}

	\subsection{Strong Convergence}
	Here we simulate the evolution of the solution on the interval $ [0,0.25] $ using a reference solution calculated with the timestep $ 10^{-6} $ and a spatial mesh with $ 64 $ frequencies. We then simulate $ 33$ initial data and $ 33$ paths for each approximation timestep. We average over all the pathwise errors to obtain the strong error of each scheme. 
	
	Figure \ref{fig:strong_convergence} shows that while the Crank-Nicolson scheme has an order of convergence close to $ 2 $, higher than the Resonance-Based scheme which shows order $ 1 $, the error is much higher for all the timesteps we tested. Hence the greater order of convergence does not compensate for the increased error in the low regularity setting. Figure \ref{fig:strong_convergence} also shows that with sufficient regularity the schemes have nearly identical order of convergence and magnitude of error.

	\subsection{Conservation Properties}
	In Figure \ref{fig:conservation_as_function_of_refinement} we study the dependence of conservation properties on the spatial mesh with the timestep $ t^{-4}$ . These plots use $ H^1 $ initial data and conservative noise. Note that the Crank-Nicolson scheme conserves mass much better than the resonance-based scheme, however, as the plots show, the resonance midpoint scheme has improved conservation properties when we increase the size of the mesh. In contrast, for the Crank-Nicolson scheme, when the mesh is too big one must decrease the timestep. This leads to divergence of the fixed point algorithm when $N=512$. This is a crucial point in the study of dispersive phenomena on long time scales because much smaller timesteps in combination with a more refined spatial mesh is computationally expensive. 
	
	In Figure \ref{fig:evo_real_H1}, using $ 10 $ initial data and $ 10 $ sample paths, we plot the sample-averaged energy and mass evolution under conservative noise. Both schemes behave comparably in this experiment.	In Figures \ref{fig:evo_real_smooth} we perform the same experiment but with smooth initial data, in which case the two schemes are exactly the same. 
	
	Finally, in Figure \ref{fig:evo_complex_smooth},  we plot the evolution of energy from smooth initial data with non-conservative (i.e. complex multiplicative) noise; again we use $ 10 $ path and data samples. Here we see that both schemes behave the same way, validating our scheme's behaviour against a traditional method. Similar experiments with other non-conservative noise show the expected behaviour for both schemes, but we do not plot them because it would be redundant.
	
	\subsection{Temporal Evolution}
	Finally, in Figure \ref{fig::evolution_plot_100}, we show a plot of the initial and end states of the evolution of the SNLSE over a long time frame. We use a spatial mesh of $ 4096 $ frequencies and a timestep of $ t^{-2} $ for a simulation on the time interval $ [0, 100] $. The cumulative error in the mass is $ 10^{-8} $, which can be improved to $ 10^{-10} $ by doubling the spatial mesh. This simulation would be unfeasible with the Crank-Nicolson scheme because one would need to use a very small timestep to obtain stability for such a refined mesh.  In Figure \ref{fig::evolution_plot_1000} we study an even larger timescale $ [0,1000] $ with the same timestep and number of frequencies. In this we see that there is a much greater dispersion effect.

\begin{figure}
	\centering
	\begin{tikzpicture}
		\begin{groupplot}[
			group style = {
				group size=2 by 1,
				horizontal sep=2cm, 
			},
			title = {Error as function of spatial refinement},
			width=0.5\textwidth,
			height=0.5\textwidth,
			grid=both,
			minor tick num=1,
			xmin = 20, xmax = 1000,
			xtick={32, 64, 128, 256, 512},
			xticklabels={32, 64, 128, 256, 512},
			major grid style={lightgray},
			minor grid style={lightgray!25},
			xlabel = {Spatial Resolution},
			legend style={
				at={(1.1,-0.3)},  
				anchor=north,
				legend columns=1,
				draw=none,
				font=\small,
			},
			]
			\nextgroupplot[
			title style={yshift=+0.2cm}, 
			ylabel={Mass Error},
			xmode = log,
			ymode = log,
			]   
			\addplot coordinates {
				(32, 1.458830e-4)
				(64, 1.724305e-5)
				(128, 2.302618e-6)
				(256, 4.334027e-7)
				(512, 1.559888e-7)
			};
			\label{plot:massRMP}
			\addlegendentry{Resonant Method}
			
			\addplot coordinates {
				(32, 1.860814e-15)
				(64, 3.849149e-15)
				(128, 2.154358e-15)
				(256, 1.678703e-14)
			};
			\label{plot:massCN}
			\addlegendentry{Crank-Nicolson}
			
			\nextgroupplot[
			title style={yshift=+0.2cm}, 
			ylabel={Energy Error},
			xlabel = {Spatial Resolution},  
			xmode = log,
			ymode = log,
			]   
			\addplot coordinates {
				(32, 4.799401e-2)
				(64, 3.035600e-2)
				(128, 2.289594e-2)
				(256, 1.500302e-2)
				(512, 1.070616e-2)
			};
			\label{plot:energyRMP}
			
			\addplot coordinates {
				(32, 6.977547e-2)
				(64, 4.234563e-2)
				(128, 3.065187e-2)
				(256, 2.127140e-2)
			};
			\label{plot:energyCN}
		\end{groupplot}
	\end{tikzpicture}
	\caption{Comparison $ L^2 $ error of mass and energy as functions of spatial refinement. Averaged over $ 33 $ initial data in $ H^1 $ and $ 33 $ paths. The noise term is multiplicative and real valued. Convergence of the fixed point method fails for Crank-Nicolson at $ N = 512 $.}
	\label{fig:conservation_as_function_of_refinement}
\end{figure}

\begin{figure}
	\centering
	\begin{tikzpicture}
		\begin{groupplot}[
			group style = {
				group size=2 by 1,
				horizontal sep=2cm, 
			},
			width=0.5\textwidth,
			height=0.5\textwidth,
			grid=both,
			minor tick num=1,
			xmin = 0, xmax = 1,
			major grid style={lightgray},
			minor grid style={lightgray!25},
			xlabel={$t$},
			legend style={
				at={(1.1,-0.2)},  
				anchor=north,
				legend columns=1,
				draw=none,
				font=\small,
			},
			]
			\nextgroupplot[
			title={Averaged Evolution of the Energy in $H^1$},
			title style={yshift=+0.2cm}, 
			ylabel={Energy},
			xmin=0, xmax=1,
			]
			\addplot table[x expr=\coordindex/9, y index=0] {out_RMP_energy_traj_H1_conservative.txt};
			\addlegendentry{Resonant Method}
			
			\addplot table[x expr=\coordindex/9, y index=0] {out_CN_energy_traj_H1_conservative.txt};
			\addlegendentry{Crank-Nicolson}
				
			\nextgroupplot[
			title={Averaged Evolution of the Mass in $H^1$},
			title style={yshift=+0.2cm}, 
			ylabel={Mass},
			ymin=1-1e-4, ymax=1+1e-4,
			xmin=0, xmax=1,
			]
			\addplot table[x expr=\coordindex/9, y index=0] {out_RMP_mass_traj_H1_conservative.txt};
			
			\addplot table[x expr=\coordindex/9, y index=0] {out_CN_mass_traj_H1_conservative.txt};
		\end{groupplot}
	\end{tikzpicture}
	\caption{Evolution of energy under influence of conservative noise with initial data in $C^\infty$.}
	\label{fig:evo_real_H1}
\end{figure}

%
%
%

\begin{figure}
	\centering
	\begin{tikzpicture}
		\begin{groupplot}[
			group style = {
				group size=2 by 1,
				horizontal sep=2cm, 
			},
			width=0.5\textwidth,
			height=0.5\textwidth,
			grid=both,
			minor tick num=1,
			xmin = 0, xmax = 1,
			major grid style={lightgray},
			minor grid style={lightgray!25},
			xlabel={$t$},
			legend style={
				at={(1.1,-0.2)},  
				anchor=north,
				legend columns=1,
				draw=none,
				font=\small,
			},
		]
		\nextgroupplot[
			title={Averaged Evolution of the Energy in $C^\infty$},
			title style={yshift=+0.2cm}, 
			ylabel={Energy},
			xmin=0, xmax=1,
			ymin=-0.01, ymax=0.05,
		]

		\addplot[solid, thick, blue, mark=o, mark repeat=2, mark phase=2] 
		table[x expr=\coordindex/9, y index=0]  {out_RMP_energy_traj_smooth.txt};
		\addlegendentry{Resonant Method }

		\addplot[red, mark=square, mark repeat=2, mark phase=0] 
		table[x expr=\coordindex/9, y index=0]  {out_CN_energy_traj_smooth.txt};
		\addlegendentry{Crank-Nicolson }
		
		\nextgroupplot[
			title={Averaged Evolution of the Mass in $C^{\infty}$},
			title style={yshift=+0.2cm}, 
			ylabel={Mass},
			ymin=0.000022, ymax=0.000024,
			xmin=0, xmax=1,
		]
		\addplot[solid, thick, blue, mark=o, mark repeat=2, mark phase=2] 
		table[x expr=\coordindex/9, y index=0]  {out_RMP_mass_traj_smooth.txt};
		
		\addplot[red, mark=square, mark repeat=2, mark phase=0] 
		table[x expr=\coordindex/9, y index=0]  {out_CN_mass_traj_smooth.txt};
	\end{groupplot}
	\end{tikzpicture}
	\caption{Evolution of energy under influence of conservative noise with initial data in $C^\infty$.}
	\label{fig:evo_real_smooth}
\end{figure}


\begin{figure}
	\centering
	\begin{tikzpicture}
		\begin{groupplot}[
			group style = {
				group size=2 by 1,
				horizontal sep=2cm, 
			},
			width=0.5\textwidth,
			height=0.5\textwidth,
			grid=both,
			minor tick num=1,
			xmin = 0, xmax = 1,
			major grid style={lightgray},
			minor grid style={lightgray!25},
			legend pos=north west,
			xlabel={$t$},
			legend style={
				at={(1.1,-0.2)}, 
				anchor=north,
				legend columns=1,
				draw=none,
				font=\small,
			},
		]
		\nextgroupplot[
			title={Averaged Evolution of the Mass in $C^{\infty}$},
			title style={yshift=+0.2cm}, 
			ylabel={Energy},
		]

		\addplot[solid, thick, blue, mark=o, mark repeat=2, mark phase=2] 
		table[x expr=\coordindex/9, y index=0]  {out_RMP_energy_traj_smooth_complex.txt};
		\addlegendentry{Resonant Method }

		\addplot[red, mark=square, mark repeat=2, mark phase=0] 
		table[x expr=\coordindex/9, y index=0]  {out_CN_energy_traj_smooth_complex.txt};
		\addlegendentry{Crank-Nicolson }
		
		\nextgroupplot[
		title={Averaged Evolution of the Mass in $C^{\infty}$},
		title style={yshift=+0.2cm}, 
		ylabel={Mass},
		]
		\addplot[solid, thick, blue, mark=o, mark repeat=2, mark phase=2] 
		table[x expr=\coordindex/9, y index=0] {out_RMP_mass_traj_smooth_complex.txt};
		\addplot[red, mark=square, mark repeat=2, mark phase=0] 
		table[x expr=\coordindex/9, y index=0] {out_CN_mass_traj_smooth_complex.txt};
		\end{groupplot}
	\end{tikzpicture}
	\caption{Evolution of the energy and mass under the influence of complex noise with initial data in $C^\infty$.}
	\label{fig:evo_complex_smooth}

\end{figure}

\begin{figure}[htbp]	
	\centering
	\begin{tikzpicture}
		\begin{groupplot}[
			group style={
				group size=2 by 1,
				horizontal sep=2cm, 
			},
			width=0.5\textwidth,
			height=0.5\textwidth,
			grid=both,
			minor tick num=1,
			xmin = 1950, xmax = 2150,
			major grid style={lightgray},
			minor grid style={lightgray!25},
			xlabel={Index},
			xtick distance=25, 
			xticklabel style={rotate=45, anchor=east, font=\small}, 
			legend style={
				at={(1.1,-0.3)},  
				anchor=north,
				legend columns=1,
				draw=none,
				font=\small,
			},
			]

			\nextgroupplot[
			title={ Magnitude },
			ylabel={Magnitude $|u|$},
			]
			\addplot[blue, opacity=0.6] table[x=index, y=init_abs] {zoomed_data_100.csv};
			\addlegendentry{Initial state}
			
			\addplot[red, opacity=0.6] table[x=index, y=end_abs] {zoomed_data_100.csv};
			\addlegendentry{Final state}

			\nextgroupplot[
			title={ Real Component of Frequencies},
			ylabel={Real part $Re(u)$},
			]
			\addplot[blue, opacity=0.6] table[x=index, y=init_real] {zoomed_data_100.csv};
			
			\addplot[red, opacity=0.6] table[x=index, y=end_real] {zoomed_data_100.csv};
			
		\end{groupplot}

	\end{tikzpicture}
	\caption{Evolution of the complex-valued vector from $t_0=0$ to $T=100$ with $ 4096 $ frequencies and timestep $ t=10^{-2} $. The plot only shows the frequencies in the range 1950 to 2150.}
	\label{fig::evolution_plot_100}
\end{figure}

\begin{figure}[htbp]	
	\centering
	\begin{tikzpicture}
		\begin{groupplot}[
			group style={
				group size=2 by 1,
				horizontal sep=2cm, 
			},
			width=0.5\textwidth,
			height=0.5\textwidth,
			grid=both,
			xmin = 1950, xmax = 2150,
			minor tick num=1,
			major grid style={lightgray},
			minor grid style={lightgray!25},
			legend pos=north west,
			xlabel={Index},
			xtick distance=25, 
			xticklabel style={rotate=45, anchor=east, font=\small}, 
			legend style={
				at={(1.1,-0.3)},  
				anchor=north,
				legend columns=1,
				draw=none,
				font=\small,
			},
			]
			
			\nextgroupplot[
			title={ Magnitude},
			ylabel={Magnitude $|u|$},
			]
			\addplot[blue, opacity=0.6] table[x=index, y=end_abs] {zoomed_data.csv};
			\addlegendentry{Initial state}
			\addplot[red, opacity=0.6] table[x=index, y=init_abs] {zoomed_data.csv};
			\addlegendentry{Final state}
			
			\nextgroupplot[
			title={ Real Component of Frequencies},
			ylabel={Real part $Re(u)$},
			]
			
			\addplot[blue, opacity=0.6] table[x=index, y=end_real] {zoomed_data.csv};
			\addplot[red, opacity=0.6] table[x=index, y=init_real] {zoomed_data.csv};
			
		\end{groupplot}
	\end{tikzpicture}
	\caption{Evolution of the complex-valued vector from $t_0=0$ to $T=1000$ with $ 4096 $ frequencies and timestep $ t=10^{-2} $. The plot only shows the frequencies in the range 1950 to 2150.}
	\label{fig::evolution_plot_1000}
\end{figure}

	\section{Conclusion}
	We have derived a class of Runge-Kutta schemes for SPDEs that reduce the regularity requirements 
	on the initial conditions compared to traditional schemes. Focusing on the example of the cubic Schr\"odinger 
	equation, we have demonstrated that, in fact, due to their implicit nature, the existence and uniqueness of 
	the schemes rely on employing a very restrictive step size, \eqref{eq::bound-on-tn}. We have shown that the 
	schemes are symplectic, and for the stochastic midpoint rule, we have proven that the local error is of 
	order $t^{3/2}$. In \cite{B.D2004}, de Bouard and Debussche remark that ``this problem [of non-uniqueness] never caused any trouble in the implementation of the method", referring to the experiments in \cite{D.D2002}. 
	Likewise, we find that our schemes are stable and convergent in practice, even though they are almost surely unstable after a long enough period of time. This mismatch between theory and practice is likely the 
	result of finite computational resources, which prevent the schemes from diverging. The resonance scheme shows better stability properties than the Crank-Nicolson suggesting the timestep is less restrictive, probably due to a better contraction constant in the stability analysis. The conclusions of our experiments are 
	that traditional methods excel at low spatial resolution and small time steps, while the resonance-based methods excel at high spatial resolution and large time steps.
	This comes at the cost of increased complexity in implementation and reduced conservation properties, which can be overcome by refining the spatial mesh -- a fair trade if 
	ones intention is to study dispersive properties.

\end{document}